\documentclass[final,12pt,a4paper]{amsart}

\title[A geometric model for the module category of a string algebra]{A geometric model for the module category of a string algebra}

\author{Karin Baur}
\email{ka.baur@me.com}
\address{Fakult\"at f\"ur Mathematik, Ruhr-Universit\"at Bochum} \address{School of Mathematics, University of Leeds}

\author{Raquel Coelho Sim\~oes}
\email{r.coelhosimoes@lancaster.ac.uk}
\address{Department of Mathematics and Statistics, University of Lancaster}

\usepackage[dvipsnames]{xcolor}
\usepackage{amssymb}
\usepackage{latexsym}
\usepackage{amsmath}
\usepackage{mathrsfs}
\usepackage[all]{xy}
\usepackage{enumerate}	
\usepackage{paralist}
\usepackage{graphicx}
\usepackage{mathabx}

\usepackage{verbatim}

\usepackage{float}   
\usepackage{rotating} 

\usepackage{hyperref}

\newcommand{\hyref}[2]{ \hyperref[#2]{#1~\ref*{#2}} }

\usepackage[notref, notcite]{showkeys} 

\theoremstyle{plain}
\newtheorem{theorem}{Theorem}[section]

\newtheorem{lemma}[theorem]{Lemma}
\newtheorem{corollary}[theorem]{Corollary}
\newtheorem{proposition}[theorem]{Proposition}

\newtheorem{introtheorem}{Theorem}

\theoremstyle{definition}
\newtheorem{remark}[theorem]{Remark}
\newtheorem{example}[theorem]{Example}
\newtheorem{definition}[theorem]{Definition} 

\newtheorem{notation}[theorem]{Notation}

\DeclareMathAlphabet{\mathpzc}{OT1}{pzc}{m}{it}

\newcommand{\sC}{\mathsf{C}}

\newcommand{\sI}{\mathsf{I}}

\newcommand{\sL}{\mathsf{L}}
\newcommand{\sM}{\mathsf{M}}
\newcommand{\sN}{\mathsf{N}}

\newcommand{\sP}{\mathsf{P}}

\newcommand{\sR}{\mathsf{R}}
\newcommand{\sS}{\mathsf{S}}

\newcommand{\sX}{\mathsf{X}}

\renewcommand{\sX}{\mathsf{X}}

\newcommand{\Fac}{\mathsf{Fac}}
\newcommand{\Sub}{\mathsf{Sub}}

\renewcommand{\geq}{\geqslant}
\renewcommand{\leq}{\leqs}
\renewcommand{\phi}{\varphi}
\renewcommand{\epsilon}{\varepsilon}

\DeclareMathOperator{\Hom}{\mathsf{Hom}}

\renewcommand{\mod}[1]{\mathsf{mod}(#1)}

\newcommand{\kk}{{\mathbf{k}}}

\newcommand{\leqs}{\leqslant}

\newcommand{\arc}[1]{\mathtt{#1}}

\entrymodifiers={+!!<0pt,\fontdimen22\textfont2>}

\newcommand{\arup}{\ar@/^/[u]}
\newcommand{\ardn}{\ar@/^/[d]}

\newcommand{\coloneqq}{\mathrel{\mathop:}=}

\usepackage{tikz,tkz-euclide}
        
\tikzset{vertex/.style={circle,fill=black,inner sep=1pt,outer sep=2pt},
         tinyvertex/.style={font=\scriptsize,minimum size=6pt},
         smallvertex/.style={inner sep=1pt, font=\small},
         >=stealth',
         leadsto/.style={-angle 90,decorate,decoration=snake,very thick},
         cut/.style={decorate,decoration=saw,very thick}}

\begin{document}

\begin{abstract}
In this paper, we give a geometric construction of string algebras and of their module categories. Our approach uses dissections of punctured Riemann surfaces with extra data at marked points, called labels. 
As an application, we give a classification of support $\tau$-tilting modules in terms of arcs in such a tiled surface. In the case when the string algebra is gentle, we recover the classification given in \cite{HZZ}.
\end{abstract}



\maketitle

{\small
\setcounter{tocdepth}{1}
\tableofcontents
}

\addtocontents{toc}{\protect{\setcounter{tocdepth}{-1}}}  
\section*{Introduction} 
\addtocontents{toc}{\protect{\setcounter{tocdepth}{1}}}   

Gentle algebras were introduced by Assem and Skowro\'nski in the context of iterated tilted algebras of type $A$ and type $\tilde{A}$ \cite{AS87}. As such, gentle algebras are one level of complexity up in the family of algebras that are morally of type $A$. They are a particularly nice class of algebras. For example, they are tame and derived-tame \cite{BM}, Gorenstein \cite{GR05} and closed under derived equivalence \cite{SZ03}. Within the class of tame algebras, the tractability of gentle algebras means they are a useful class of algebras to serve as test cases for more general conjectures in representation theory. In recent years, there has been renewed interest in gentle algebras due to their ubiquity, having connections with symplectic geometry via partially wrapped Fukaya categories \cite{HKK17, LP}, algebraic geometry and theoretical physics via dimer models \cite{Bro}, Lie theory \cite{HK01} and cluster theory.  

Triangulations of Riemann surfaces give a geometric framework to study cluster theory \cite{FST}. This framework provides geometric models giving detailed information regarding the representation theory of associated (Jacobian) algebras and associated (cluster) categories (see e.g.~\cite{ABCP, BM08, BZ11, CCS06, LF09}). Since then, geometric models have been developed in many contexts. For example, geometric models have been used to describe not only cluster categories but also module categories of gentle algebras \cite{BCS21}, derived categories of gentle algebras and their connection with Fukaya categories of surfaces \cite{LP, OPS18}, module categories and derived categories of skew-gentle algebras \cite{HZZ, LFSV}, a generalisation of gentle algebras which are morally of type $D$, and very recently module categories of semilinear locally gentle algebras \cite{BBTJW}. 

String algebras are the natural next step up in complexity in the family of algebras of ``type $A$''. In particular, all gentle algebras are string algebras. Like gentle algebras, the building blocks of the representation theory over string algebras, such as indecomposable modules, irreducible morphisms and Auslander-Reiten sequences, are well understood \cite{BR87, CB89, Krause}. However, there are notable differences. For example, string algebras need not be Gorenstein and while they are tame they need not be derived-tame. It is natural to ask therefore whether module categories of string algebras can be described using geometric models and whether these geometric models can be used to shed light on the less well studied representation theory of string algebras and to further develop connections with other areas of mathematics. 

Gentle algebras can be realised in terms of so-called {\it tiling algebras} associated to dissections of unpunctured surfaces with boundary. This realisation can be extended to locally gentle algebras, where we drop the finite-dimensional assumption. These can be associated to dissections of surfaces which may be punctured and may not have boundary \cite{PPP19}. By viewing a string algebra as a quotient of a locally gentle algebra, we realise string algebras as so-called {\it labelled tiling algebras} associated to dissections of punctured surfaces, possibly with no boundary, together with some extra data, in terms of labels at marked points.

\begin{introtheorem}[Theorem~\ref{thm:model-algebra}]\label{thmA}
Let $A$ be a finite dimensional monomial algebra. Then $A$ is a string algebra if and only if $A$ is a labelled tiling algebra of a marked punctured surface.
\end{introtheorem}

We are then able to describe the building blocks of the representation theory of a string algebra combinatorially using an associated labelled tiling (or dissection) of a surface. 

\begin{introtheorem}[Theorems~\ref{thm:permissiblearcsstrings} and~\ref{thm:pivot-move-irreducible}]\label{thmB}
Let $A=A(\sS)$ be a labelled tiling algebra with associated surface $\sS$. There are explicit bijections between:
\begin{compactenum}
\item the equivalence classes of permissible arcs in $\sS$ and the strings of $A$; 
\item isotopy classes of certain permissible closed curves in $\sS$ and powers of bands; 
\item certain moves between permissible arcs and irreducible morphisms of string modules.
\end{compactenum}
\end{introtheorem} 

In the case where there are no labels, the associated tiling algebra is locally gentle. We recover the gentle algebras case when there are no labels and the endpoints of the dissection lie on the boundary of $\sS$. 

Certain permissible arcs, called $\sR$-arcs (see Definition~\ref{def:M-R-arcs}), permit us to give a combinatorial description of all morphisms between string modules via crossings and intersections (see Theorem~\ref{thm:R-arcs-morphisms}, and \cite{Chang} for the gentle algebras case). In upcoming \cite{BCSGS24}, $\sR$-arcs are used to study maximal almost rigid objects over gentle algebras. These are new representation theoretic objects defined in \cite{BGMS}, for algebras of type $A$, in connection with Cambrian lattices. 

Geometric models not only play an important role in relating different areas of mathematics, they are also useful to study quite often long-standing questions in representation theory. For example, geometric models for gentle algebras have been applied to the study of:
\begin{itemize}
    \item $d$-finite representation algebras in higher homological algebra \cite{HJS};
    \item derived invariants for algebras \cite{AG, ALFP, APS, DR, LP, Opper, OZ};
    \item length hearts which are crucial to the construction of Bridgeland stability conditions \cite{Chang};
    \item silting theory \cite{APS, CJS, CS23}, in particular in providing counterexamples for silting completion \cite{LZ,JSW}.
\end{itemize}

Geometric models have also been useful in $\tau$-tilting theory, introduced in \cite{AIR} as a ``completion'' of classical tilting theory from the point of view of mutation. The geometric models for gentle and skew gentle algebras were used to classify support $\tau$-tilting modules over these algebras (\cite{HZZ}, see also \cite{BDMTY, PPP18}) and to study the connectedness of the mutation graph \cite{HZZ22}. 

In this paper, we classify support $\tau$-tilting modules over string algebras, using a different set of permissible arcs, called $\sM$-arcs (see Definition~\ref{def:M-R-arcs}). In the gentle algebras case, considered in \cite{HZZ}, support $\tau$-tilting modules can be described in terms of non-crossing arcs.  For string algebras in general we have to allow some crossings, but these are easily identifiable by a local condition with respect to labels. This is what we mean by ``good crossings'' below (see Definition~\ref{good-crossing} for the precise definition).

\begin{introtheorem}[Theorem~\ref{thm:tau-tilt}]\label{thmC}
    There is a bijection between the set of support $\tau$-tilting modules and the set of maximal collections of generalised permissible arcs with good crossings. 
\end{introtheorem}

\subsection*{Acknowledgments} 
The authors would like to thank the Isaac Newton Institute for Mathematical Sciences, Cambridge, for support and hospitality during the programme {\it Cluster algebras and representation theory}, where work on this paper was undertaken. This work was supported by EPSRC grant EP/R014604/1. This project has also been supported by the EPSRC Programme Grant EP/W007509/1 and by the Royal Society Wolfson Award RSWF/R1/180004 of the first author and by European Union’s Horizon 2020 research and innovation programme through the Marie Sk\l odowska-Curie Individual Fellowship grant 838706 of the second author.

\section{Background on string algebras}\label{sec:background}

Let $\kk$ be an algebraically closed field and $Q = (Q_0, Q_1)$ be a quiver, where $Q_0$ is the set of vertices and $Q_1$ is the set of arrows. Let $I$ be an admissible ideal of $\kk Q$, and $A$ be the bound quiver algebra $\kk Q/I$. 

Given $a \in Q_1$, we denote its source by $s(a)$ and its target by $t(a)$. We will read paths in $Q$ from left to right, and $\mod A$ will denote the category of right $A$-modules. 

The finite dimensional algebra $A$ is said to be {\it monomial} if $I$ is generated by paths of length at least two. We refer to a path in $I$ of length $k$ a {\it $k$-relation}. 

\begin{definition}\cite{AS87}\label{def:string} 
A finite dimensional algebra $A$ is a {\it string} algebra if it admits a presentation $A = \kk Q/I$ satisfying the following conditions: 
\begin{compactenum}[(S1)]
\item Each vertex of $Q$ is the source of at most two arrows and the target of at most two arrows.
\item For each arrow $a$ in $Q$, there is at most one arrow $b$ in $Q$ such that $ab \not\in I$, and there is at most one arrow $c$ such that $ca \not\in I$.  
\item $I$ is generated by paths of length at least two, i.e. $A$ is monomial.
\end{compactenum}
\end{definition}

\begin{definition}\label{def:gentle}
A finite dimensional algebra $A=\kk Q/I$ is {\it gentle} if is it a string algebra satisfying the following additional conditions:
\begin{itemize}
\item[(G1)]  For each arrow $a$ in $Q$, there is at most one arrow $b'$ such that $ab'\in I$ and there is at most one arrow $c'$ such that $c'a\in I$. 
\item[(G2)] $I$ is generated by $2$-relations. 
\end{itemize} 
\end{definition}

A vertex $v$ satisfying (S1) and for which the paths of length 2 going through $v$ satisfy (S2) and (G1) is called a {\it gentle vertex}. 

\begin{definition}\label{def:locally-gentle}
A {\it locally gentle algebra} is an algebra which admits a presentation $\kk Q/I$ satisfying the conditions in Definition~\ref{def:gentle} but which is not necessarily finite dimensional, i.e. $I$ is not necessarily admissible.
\end{definition}

\begin{example}\label{ex:string-gentle}
Consider the quiver $Q$
\[
\xymatrix{
4\ar[r]^{d}  & 1\ar[r]^{e}\ar[d]^{a}  & 5 \\
3\ar[u]^{c} & 2 \ar[l]^{b}
}
\] 
and the following ideals: 
\[
I_1=\left<de, da, abc \right>, 
\ 
I_2=\left<de \right>, 
I_3=\left<da\right>
\]
Then $\kk Q/I_1$ is a string algebra which is not gentle nor locally gentle. The algebra 
$\kk Q/I_2$ is locally gentle but not gentle and 
$\kk Q/I_3$ is gentle. 
See also Example~\ref{ex:two-surfaces-one-algebra} where we consider a different ideal. 
\end{example}

\subsection{Strings and bands} 
Let $A=\kk Q/I$ be a string algebra. Given $a \in Q_1$, we define a formal inverse $a^{-1}$ such that $s(a^{-1}) = t(a)$ and $t(a^{-1}) = s(a)$. We denote by $Q_1^{-1}$ the set of formal inverses of arrows of $Q$. 

A {\it walk} is a sequence $w= a_1 \ldots a_r$, with $a_i \in Q_1^{-1} \cup Q_1$ for each $i$ and $t(a_i) = s(a_{i+1})$. Each of the $a_i$'s is called a {\it letter of $w$}. Given a walk $w=a_1 \ldots a_r$, a {\it subwalk of $w$} is a sequence of the form $a_i a_{i+1} \ldots a_j$ for some $1 \leq i \leq j \leq r$, and the {\it inverse of $w$} is the walk $a_r^{-1} \ldots a_1^{-1}$. 

A {\it string} is a walk which is {\it reduced}, i.e. it has no subwalks of the form $aa^{-1}$ or $a^{-1}a$ with $a \in Q_1$, and avoids relations, i.e. it has no subwalks $v$ with $v \in I$ or $v^{-1} \in I$. If all letters of a string lie in $Q_1$ (resp. $Q^{-1}_1$), we say the string is {\it direct} (resp. {\it inverse}). 

Each vertex $i \in Q_0$ defines two {\it trivial strings}, denoted by $1_i^+$ (or simply $1_i$) and its inverse $1_i^{-}$. For technical reasons, we also need to define the {\it zero string}, which will be denoted by $w=0$. 

Given two non-trivial (and non-zero) strings $w=a_1 \ldots a_r$ and $v=b_1 \ldots b_s$, the concatenation $wv = a_1 \ldots a_r b_1 \ldots b_s$ is well defined if $t(a_r) = s(b_1)$. Moreover, $wv$ is a string if and only if $b_1 \neq a_r^{-1}$ and there is no relation of the form $a_i \ldots a_r b_1 \ldots b_j$, for some $1\leq i \leq r$ and $1 \leq j \leq s$. 

Now suppose one of the strings, say $w$, is trivial. Write $w= 1_x^\gamma$, where $\gamma \in \{+,-\}$ and $x \in Q_0$. We would like there to be at most one choice of arrow $a$ (resp. $b$) with $s(a) = x$ (resp. $t(b) = x$) for which $1_x^\gamma a$ (resp. $1_x^\gamma b^{-1}$) is a string. However, there could be two arrows starting (resp. ending) at $x$. For this reason, we need the following technical definition.

\begin{definition}\label{def:sign-functions}
    Given a string algebra $A= \kk Q/I$, we define two {\it sign functions} $\sigma, \epsilon \colon Q_1 \rightarrow \{-1,1\}$ satisfying the following conditions:
\begin{itemize}
    \item If $b_1, b_2 \in Q_1$ are distinct arrows with the same start, then $\sigma(b_1) \neq \sigma (b_2)$.
    \item If $a_1, a_2 \in Q_1$ are distinct arrows with the same target, then $\epsilon (a_1) \neq \epsilon (a_2)$.
    \item If $a, b \in Q_1$ with $s(b) = t(a)$ and $ab \not\in I$, then $\sigma(b) = -\epsilon (a)$.
\end{itemize}

We can extend the domain of these functions to the set of all strings, as follows. 

If $b\in Q_1$, we define $\sigma (b^{-1}) \coloneqq \epsilon (b)$ and $\epsilon (b^{-1}) \coloneqq \sigma (b)$. 

If $w = a_1 \ldots a_n$ is a non-trivial string (i.e. $n\geq 1$), then $\sigma (w) \coloneqq \sigma (a_1)$ and $\epsilon (w) \coloneqq \epsilon (a_n)$. 

Given $v \in Q_0$, we define $\sigma (1_v^\pm) \coloneqq \mp 1$ and $\epsilon (1_v^\pm) \coloneqq \pm 1$.
\end{definition}

We are now ready to extend the definition of concatenation of strings to involve trivial strings as well. Let $v=1_x^\gamma$, where $x \in Q_0$ and $\gamma \in \{1,-1\}$, and let $w$ be a non-trivial string. Then $vw$ is a string if and only if $s(w) = x$ and $\sigma(w) = -\epsilon (v) = -\epsilon(1_x^{\gamma}) = -\gamma$. Similarly, $wv$ is a string if and only if $t(w) =x$ and $\epsilon (w) = - \sigma (v) = \gamma$. 

\begin{remark}
The reader can find an algorithm for choosing sign functions in~\cite[p.~158]{BR87}. We note that if $A$ is (locally) gentle, we always have a choice of $\sigma, \epsilon$ so that given two strings $v, w$, the concatenation $vw$ is defined and is a string if and only if $t(v) = s(w)$ and $\sigma (w) = - \epsilon (v)$. However, this is not the case for string algebras in general. For example, consider the quiver
\[
\xymatrix@R.9pt{
& & & 4\\
1 \ar[r]^a & 2 \ar[r]^b & 3 \ar[ur]^c \ar[dr]_d & \\
& & & 5
}
\]
bound by the relations $abd$ and $bc$. For any choice of $\sigma$ and $\epsilon$, we must have $\sigma (d) = - \epsilon (b)$ because $bd$ is a string. However, it follows that $\sigma (d) = - \epsilon (ab)$, but $abd$ is not a string. 
\end{remark}

Every string $w$ of $A$ defines an indecomposable module $M(w)$, called {\it string module}. If we realise the string $w$ as a graph, the string module $M(w)$ is given by assigning a copy of the field $\kk$ to each of the vertices of the graph $w$ and a copy of the identity map to each of the arrows in graph. Given two strings $w,v$, we have $M(w) \simeq M(v)$ if and only if $w=v$ or $w=v^{-1}$. The module corresponding to the zero string is the zero module.

A string $b=b_1 \ldots b_r$ which is cyclic, i.e. $t(b_r) = s(b_1)$, and such that any power $b^k$ of $b$ is a string but $b$ itself is not a proper power of any string, is called a {\it band}. 

Each band $b = b_1 \ldots b_r$ (up to cyclic permutation and inversion) defines a family of indecomposable modules $M(b,n,\varphi)$, called {\it band modules}, where $n\in \mathbb{N}$ and $\varphi$ is an irreducible automorphism of $\kk^n$. The band module $M(b,n, \varphi)$ is defined by assigning a copy of $\kk^n$ to each vertex of $b$, the identity map to each arrow of $b$ apart from $b_r$ which is assigned the automorphism $\varphi$.

\begin{theorem}\cite{BR87}
    Every indecomposable module over a string algebra is either a string module or a band module.
\end{theorem}

%
\subsection{Description of morphisms between string modules}\label{sec:morphisms-strings}

We recall here how morphisms can be described in terms of string combinatorics, using the language of~\cite{BCS21-corr}. 

Let $w$ be a string. 
A {\it factor string decomposition} of $w$ is a decomposition of the form $w=w_\ell w_e w_r$, such that either $w_\ell$ is a trivial string or the last letter of $w_\ell$ lies in $Q_1^{-1}$ and either $w_r$ is trivial or the first letter of $w_r$ lies in $Q_1$. In this case, $w_e$ is called a {\it factor substring} of $w$. 

Similarly, a decomposition of $w$ of the form $w=w_\ell  w_e  w_r$, where either $w_\ell$ is a trivial string or the last letter of $w_\ell$ lies in $Q_1$  and either $w_r$ is trivial or the first letter of $w_r$ lies in $Q_1^{-1}$, is called a {\it substring decomposition}. In this case, $w_e$ is called an {\it image substring} of $w$. 

The set of equivalence classes of factor string (resp. substring) decompositions of $w$ is denoted by $\Fac (w)$ (resp. $\Sub(w)$).  

Let $v, w$ be two strings. We call a pair $(v_\ell v_e v_r, w_\ell w_e w_r)$ in $\Fac (v) \times \Sub (w)$ an {\it admissible pair} if $w_e = v_e$ or $w_e = v_e^{-1}$. We will also sometimes call the string $w_e$ an {\it overlap between $w$ and $v$}. 

Homomorphisms between string modules arise from admissible pairs as has been shown by Crawley-Boevey:  

\begin{theorem}\cite{CB89}\label{thm:homs}
Let v, w be strings and $M(v), M(w)$ be the corresponding string modules.
Then $\dim_\kk \, \Hom_A (M(v), M(w))$ is given by the number of  admissible pairs in $\Fac (v) \times \Sub (w)$. 
\end{theorem} 

\subsection{Irreducible morphisms and the Auslander-Reiten translate}
In order to describe irreducible morphisms between string modules, we need the definition of hook and cohook. 

\begin{definition}
    Let $A=\kk Q/I$ be a string algebra and let $a \in Q_1$.
    \begin{enumerate}
        \item Let $p_a$ denote the direct string with $s(p_a) = s(a)$ such that it does not start with $a$ and it is {\it right maximal}, i.e. there is no $b \in Q_1$ with $p_a b$ a string. We call the string $a^{-1} p_a$ the {\it hook of $a$}. If $s(a)$ has out-valency one, then $p_a$ is a trivial string at $s(a)$.
        \item Let $q_a$ denote the direct string with $t(q_a) = t(a)$ and such that it does not end with $a$ and it is {\it left maximal}, i.e. there is no $b \in Q_1$ with $b q_a$ a string. We call the string $a q_a^{-1}$ the {\it cohook of $a$}. If $t(a)$ has in-valency one, then $q_a$ is a trivial string at $t(a)$. 
    \end{enumerate}
\end{definition}

Let $w$ be a string. We define a new string $f_s(w)$ as follows: 
\begin{itemize}
    \item If there is $a \in Q_1$ such that $wa^{-1}$ is a string, then $f_s(w)$ is obtained by adding the hook of $a$ on the right, i.e. $f_s(w) \coloneqq wa^{-1} p_a$. 
    \item If there is no $a \in Q_1$ such that $wa^{-1}$ is a string and $w$ is not inverse, then $f_s(w)$ is obtained from $w$ by removing a cohook on the right. In other words, $w$ must be of the form $w_\ell a q_a^{-1}$, for some $a\in Q_1$, in which case $f_s(w) \coloneqq w_\ell$.
    \item If there is no $a \in Q_1$ such that $wa^{-1}$ is a string and $w$ is an inverse string, then $f_s(w) \coloneqq 0$. 
\end{itemize}

Similarly, we define $f_t(w)$ to be of the form $p_a^{-1}a w$, if there is $a\in Q_1$ with $aw$ a string (adding a hook on the left), or of the form $w_r$, where $w=q_a a^{-1}w_r$ (removing a cohook on the left), or $f_t(w) = 0$ if $w$ is direct and left maximal. 

\begin{remark}\label{rmk:not-commutative}
    Let $A$ be a string algebra and $w$ be a string over $A$. It is not necessarily the case that $f_t(f_s(w)) = f_s(f_t(w))$. For example, consider the quiver 
    \[\xymatrix{1 \ar[r]^a & 2 \ar[r]^b & 3 \ar[r]^c & 4 & 5 \ar[l]_d & 6 \ar[l]_e & 7 \ar[l]_f}
    \]
    bound by the relations $abc$ and $fed$. Consider the string $w=bcd^{-1}e^{-1}$. We have $f_s(w) = b, f_t(f_s(w)) = ab, f_t(w) = e^{-1}$ and $f_s(f_t(w)) = e^{-1} f^{-1}$. 
\end{remark}

\begin{lemma}\label{lem:non-injective}
    Let $A$ be a string algebra and $w$ a string over $A$ such that $M(w)$ is not an injective module. Then at least one of $f_s(w)$ or $f_t(w)$ is non-zero. Moreover, if both $f_s(w)$ and $f_t(w)$ are non-zero strings then $f_s(f_t(w)) = f_t(f_s(w))$.
\end{lemma}
\begin{proof}
    Since $M(w)$ is not injective, we can assume without loss of generality that there is $a\in Q_1$ such that $aw$ is a string. In particular $f_t(w) \neq 0$. Now assume $f_s(w)$ is also non-zero.

    {\bf Case 1.} There is an arrow $b$ such that $awb^{-1}$ is a string.

    In this case, we have $f_t(w) = p_a^{-1}aw$ and $f_s(f_t(w))=p_a a^{-1}wb^{-1}p_b$. On the other hand, $f_s(w) = wb^{-1} p_b$ and so $f_t(f_s (w)) = p_a a^{-1}wb^{-1}p_b = f_s(f_t(w))$. 

    {\bf Case 2.} There is an arrow $c$ such that $w=w_\ell c q_c^{-1}$. 

    Then $f_t(w) = p_a^{-1}aw$ and $f_s(f_t(w))= p_a^{-1}a w_\ell$. On the other hand, $f_s(w) = w_\ell$ and $f_t(f_s(w)) = p_a^{-1}aw_\ell$, and we are done.

The only remaining case is when $w=q_a^{-1}$. However, we would have $f_s(w) = 0$, contradicting the hypothesis. 
\end{proof}

For a string module $M(w)$, the irreducible morphisms starting at $M(w)$  and its Auslander-Reiten translate can be described in terms of the strings $f_s(w)$ and $f_t(w)$ as follows.

\begin{theorem}\label{thm:irreducible}\cite{BR87}
 Let $A$ be a string algebra and $w$ a string. 
 \begin{enumerate}
     \item There are at most two irreducible morphisms starting at $M(w)$, with possible targets $M(f_s(w))$ and $M(f_t(w))$. 
     \item If $M(w)$ is non-injective, suppose without loss of generality that $f_s(w) \neq 0$. Then the Auslander-Reiten sequence starting at $M(w)$ is of the form
     \[
     \xymatrix{0 \ar[r] & M(w) \ar[r] & M(f_s(w)) \oplus M(f_t(w)) \ar[r] & M(f_t(f_s(w)) \ar[r] & 0.}
     \]
     In particular, $\tau^{-1} (M(w)) = M(f_t(f_s(w)))$. 
 \end{enumerate}
\end{theorem}

The remaining irreducible morphisms are between band modules of the form $M(b, n, \varphi)$ and $M(b, n+1, \varphi)$. Moreover, the Auslander-Reiten translate acts as the identity morphism on band modules, since these modules lie in homogeneous tubes. 

\section{Labelled tiling algebras}\label{sec:set-up}


%
\subsection{Tiled surfaces and labels}\label{ssec:tilings}

We now recall the geometric set-up appearing in the construction of locally gentle algebras and introduce labels on tiled surfaces as a means towards string algebras. 

\begin{definition}
    A {\it marked surface} $(\sS,\sN)$ is given by a 
    connected oriented Riemann surface $\sS$ with a finite set $\sN$ of marked points on the (possibly empty) boundary  and in the interior of $\sS$. We impose that each boundary component has marked points. We  write $\sN_b$ for $\sN\cap \partial\,\sS$ and $\sN^{\circ}=\sN\setminus\sN_b$ for the marked points in the interior. The elements of $\sN^{\circ}$ are called punctures. We will often simply write $\sS$ for the tuple $(\sS,\sN)$. 
\end{definition}

The segments between two neighboured marked points on a boundary component are called the {\it boundary segments} of $(\sS,\sN)$. 

\begin{definition}\label{def:arc} 
Let $\sS=(\sS,\sN)$ be a marked surface. 
\begin{enumerate}
    \item 
    An {\it arc} $\gamma$ in $\sS$ is a curve in 
$\sS$ such that its only intersections with $\sN$ are the endpoints of $\gamma$. Curves are taken up to isotopy fixing endpoints. 
    \item 
    Two arcs $\gamma_1$ and $\gamma_2$ are {\it intersecting} if they have a common endpoint. 
    \item 
    Two arcs $\gamma_1$, $\gamma_2$ are {\it non-crossing} if there exist representatives in their isotopy classes which do not meet anywhere except possibly at  endpoints. 
    The arcs are called {\em crossing} otherwise. 
    \item 
    A collection of arcs is {\it non-crossing} if none of the arcs has self-crossings in the interior of $\sS$ and the arcs are pairwise non-crossing.
\end{enumerate}
\end{definition}

\begin{definition}\label{def:tiled}
A non-crossing collection $\sP$ of arcs in $\sS$ is called a 
{\it dissection} (also known as {\it partial triangulation} or {\it tiling}) of $\sS$ if the following holds: 
\begin{itemize}
    \item no arc of $\sP$ is isotopic to a boundary segment or to a point. 
    \item 
    The arcs of $\sP$ dissect the surface into two types of tiles described as follows:
    \begin{enumerate}
    \item[(i)] polygonal tiles of size $\ge 3$ where exactly two consecutive sides are boundary segments; 
    \item[(ii)] polygonal tiles of size $\ge 1$ with no boundary segment and with a marked point in the interior.
\end{enumerate}   
\end{itemize} 
If $\sP$ is a dissection of $\sS$, the pair $(\sS,\sP)$ is also called 
a {\it tiled surface}. 
\end{definition} 

\begin{figure}[ht!]
    \centering
    \includegraphics[scale=.6]{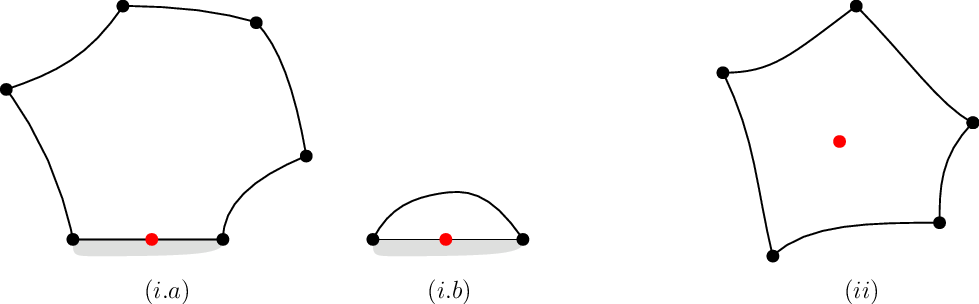}
    \caption{The  types of tiles in a dissection. Figure $(i.b)$ is the special case of a tile with boundary segments and of size $3$.}
\label{fig:tiles-red}
\end{figure}

\begin{figure}[ht!]
 \centering
\includegraphics[scale=.7]{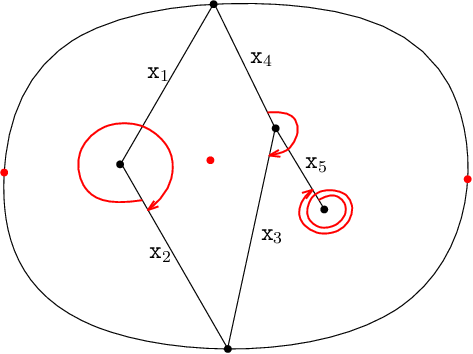}
\caption{Labeled tiling of a surface.} 
\label{fig:label-angle-red}
\end{figure}

Figure~\ref{fig:tiles-red} gives an illustration of the types of tiles. We note that tiles are allowed to be self-folded. For example, in Figure~\ref{fig:label-angle-red} the tile on the right is an hexagon with the edge $5$ glued to itself.

Note that each tile of the dissected surface contains exactly one marked point in $\sN$ which is not incident to any arc in $\sP$. If the tile has a boundary segment, then this marked point lies in it, otherwise the marked point is a puncture. As illustrated in Figure~\ref{fig:tiles-red}, we will draw these ``isolated'' marked points as red dots, and as such we will often call them {\it red points}. 

\begin{notation}
  A dissection $\sP$ of $\sS$, defines a partition $\sM \cup \sR$ of the set of marked points $\sN$ where $\sR$ denotes the set of red points, i.e. the marked points not incident with any of the arcs in $\sP$, and $\sM$ denotes the set of the remaining points, which are drawn in the figures as black dots. 
\end{notation}

We note that in~\cite{BCS21}, unmarked boundary components were allowed in certain tiles of type (ii). Here, we replace these unmarked boundary components with punctures in $\sR$. 

Let $p \in \sM$. If $p$ lies in the boundary of $\sS$, let $p'$ and $p''$ be two points in the same boundary component such that $p', p'' \not\in \sN$ and $p$ is the only marked point in $\sN$ lying in the boundary segment between $p'$ and $p''$. Consider a curve $\delta$ isotopic to this boundary segment but such that the only intersection with the boundary of $\sS$ is at its endpoints. If $p$ is a puncture, consider a closed simple curve $\delta$ in the interior of $\sS$ around $p$ and with no intersections with $\sN$. 

The {\it complete fan at $p$} is defined to be the sequence of all arcs of $\sP$ that $\delta$ crosses in the clockwise order. A {\it fan at $p$} is any subsequence of consecutive arcs of the complete fan at $p$. 

\begin{remark}
    Let $p \in \sM$.
    \begin{enumerate}
        \item If an arc in $\sP$ is a loop, then it lies in exactly one complete fan. All the remaining arcs lie in two distinct complete fans. 
        \item The complete fan at $p \in \sN^\circ \cap \sM$ is an infinite periodic sequence, with period $\geq 1$. 
    \end{enumerate}
\end{remark}

\begin{figure}[ht!]
    \centering
\includegraphics[scale=0.8]{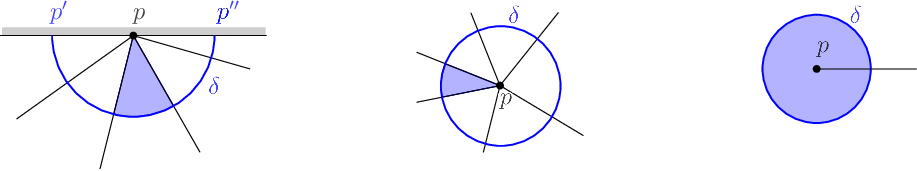}
    \caption{The shaded triangle is an angle at $p$. A puncture $p \in \sM$ can have a unique angle which is a self-folded triangle.}
    \label{fig:angle}
\end{figure}

Any two consecutive crossings of the curve $\delta$ with arcs in $\sP$ define a triangle whose vertices are these two crossings and the marked point $p$; see Figure~\ref{fig:angle}. This triangle is called an {\it angle at $p$} (also called {\it angle of $\sP$}), and we say the curve $\delta$ {\it cuts this angle at $p$}. 

We will now associate an algebra to a tiled surface. We note that the subset $\sR$ of the marked points $\sN$ of the surface does not play a role in the definition of this algebra. The subset $\sR$ will be however useful in Section~\ref{sec:geom-model}, when describing the module category of this algebra via the surface. In the context of gentle algebras, this definition recovers that in~\cite{BCS21, OPS18} and in the context of locally gentle algebras, which includes the case where the surface has no boundary, this definition recovers that in~\cite{PPP19}. 

\begin{definition}\label{def:quiver}
Let $(\sS,\sP)$ be a tiled surface. We define $Q=Q(\sS,\sP)$ to be the quiver defined as follows: 

\begin{compactenum}
\item The vertices in $Q$ are in bijection with the arcs in $\sP$.
\item The arrows in $Q$ are in bijection with the angles at the marked points in $\sM$. In other words, there is an arrow $a\colon x_i \to x_j$ in $Q$ if and only if the 
corresponding arcs $\arc{x}_i$ and $\arc{x}_j$ share an endpoint $p_a \in \sM$ and $\arc{x}_j$ is the immediate successor of $\arc{x}_i$ in the complete fan at $p_a$. 
\end{compactenum}
\end{definition}

\begin{definition}\label{def:label}
Let $(\sS,\sP)$ be a tiled surface. 
\begin{enumerate}
\item A {\it label} at a point $p\in \sM$ is a finite fan at $p$ of length $\geq 3$. 
 \item Given two labels $\ell$ and $\ell'$, if $\ell$ is a subsequence of $\ell'$, then we say the label $\ell'$ is {\it redundant}. 
\end{enumerate}
\end{definition}

\begin{remark}
    Each label at a point $p \in \sM$ defines a linearly oriented path in the quiver $Q$. If $p$ is a puncture, this path may be of the form $a^k$, where $a \in Q_1$ is such that $s(a) = t(a)$; see Figure~\ref{fig:label-angle-red} for an example.
\end{remark}

\begin{notation}
We will draw a label $\ell$ as a red curve $\delta$ around the corresponding marked point with an orientation such that $\delta (0)$ is a point in the first arc of $\ell$, $\delta (1)$ is a point in the last arc of $\ell$ and if $\ell$ contains an arc $\arc{x}$ multiple times, the corresponding intersections of $\delta$ with $\arc{x}$ are at distinct points. This is in order to make it clear what the label is from $\delta$, especially in the case when the label starts and ends at the same arc. See Example~\ref{ex:labeled} and Figure~\ref{fig:label-angle-red}. 
\end{notation}

\begin{definition}\label{def:ideal}
Let $(\sS,\sP)$ be a tiled surface and let $\sL$ be a finite set of labels. 
Let $Q=Q(\sS,\sP)$ and consider the path algebra $\kk Q$. 
We define $I_{\sP,\sL}$ to be the ideal of $\kk Q$ generated by the following paths: 
\begin{compactenum}
\item $ab$ of length 2 such that $p_a \neq p_b$ or $p_a = p_b$ and $t(a) = s(b)$ corresponds to a loop arc. 
\item linearly oriented paths in $Q$ corresponding to a label in $\sL$.
\end{compactenum}
\end{definition} 

The relations of type (1) arise from the tiles in $(\sS, \sP)$, as certain compositions of two arrows within a tile. These are all the relations considered in the construction of (locally) gentle algebras via surfaces. For string algebras in general, we need also to consider labels and corresponding relations of type (2).

\begin{definition}\label{def:tiling-alg}
Let $\sS=(\sS,\sP)$ be a tiled surface and let $Q=Q(\sS,\sP)$. 
Let $\sL$ be a finite set of labels on $\sS$ such that for every $q\in \sN^\circ \cap \sM$, $\sL$ contains a label at $q$. 
The {\it labelled tiling algebra} $A_{\sP,\sL}$ associated to the data $(\sS, \sP, \sL)$ is the bound quiver algebra $A_{\sP, \sL} = \kk Q/I_{\sP,\sL}$.
\end{definition} 

Note that the reason to impose labels at all punctures in $\sM$ is to ensure that the algebras are finite dimensional. One could drop this and work with infinite dimensional string algebras instead. 

From now on we will assume for simplicity that there are no redundant labels in $\sL$. The following example illustrates these concepts. 

\begin{figure}[ht!]
\centering    
\includegraphics[scale=.9]{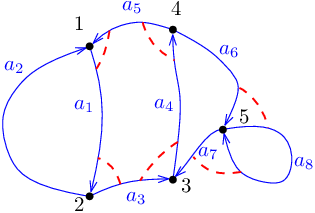}
\caption{The quiver of the tiling from Figure~\ref{fig:label-angle-red}. The type (1) relations are indicated by dashed (red) lines.}
\label{fig:quiver-tiling-red}
\end{figure}

\begin{example}\label{ex:labeled}
Consider the tiling $\sP$ of a disk $\sS$ with four punctures in Figure~\ref{fig:label-angle-red}. 
The quiver of this tiling is in Figure~\ref{fig:quiver-tiling-red}, where vertex $i$ corresponds to the arc $\arc{x}_i$, for each $1 \leq i \leq 5$.

Consider the labels given by the sequences $\{\arc{x}_2, \arc{x}_1, \arc{x}_2\}$, $\{\arc{x}_4, \arc{x}_5, \arc{x}_3\}$ and $\{\arc{x}_5, \arc{x}_5,\arc{x}_5\}$. 
The ideal $I_{\sP,\sL}$ is generated by the 2-relations of type (1)
\[
\{a_1a_3,a_3a_4, a_4a_5,a_5a_1,a_6a_8, a_8a_7\},
\]
together with the type (2) relations:
$\{a_2a_1, a_6a_7,a_8^2\}$. 
\end{example}

\section{String algebras as  labelled tiling algebras}\label{sec:labelled-tiling}

In this section we will show that the class of labelled tiling algebras coincide with the class of string algebras. 
We do this by realising a string algebra as the quotient of a locally gentle algebra, by removing all relations of length $\geq 3$ and a choice of 2-relations so that every vertex becomes gentle. The relations removed will correspond to the labels. 

\begin{theorem}\label{thm:model-algebra}
Let $A = \kk Q/I$ be a finite dimensional monomial algebra. The following are equivalent. 
\begin{enumerate}
    \item 
    $A$ is a string algebra.
    \item
    $A=(\kk Q/J)/(I/J)$, where $\kk Q/J$ is a locally gentle algebra. 
    \item
    $A$ is a labelled tiling algebra of a marked surface. 
\end{enumerate}
\end{theorem}

\begin{proof}
$(1)\Longleftrightarrow (2)$: 
Suppose $A = \kk Q/I$ is a string algebra. Pick a choice of two sign functions $\sigma, \epsilon$ on the set of strings which satisfy the conditions in Definition~\ref{def:sign-functions}. Now let $U$ be a set of generators of $I$ and let $V$ be the subset of $U$ given by all paths of length $\ge 3$ and all paths  $ab$ of length $2$ such that $\sigma(b) = -\epsilon(a)$. Let $J = \langle U \setminus V \rangle$.

Clearly we have $A = (\kk Q/J)/(I/J)$. Moreover, the algebra $B\coloneqq \kk Q/J$ is locally gentle. Indeed, if there is a vertex $x$ in $Q_0$ which is not gentle in $B$, then we can assume without loss of generality, that there are two relations in $B$ of the form $a_1 b$ and $a_2 b$, with $s(b) = x$. But then, for the choice of sign functions above, we would have $\epsilon (a_1) \neq \epsilon (a_2)$, and so $\sigma (b) = - \epsilon (a_i)$, for some $i\in \{1,2\}$. This contradicts the fact that $a_i b$ is a relation in $B$.  The converse is clear.

$(2)\Longleftrightarrow (3)$: 
This correspondence follows from the bijection between the set of isomorphism classes of locally gentle algebras and the set of homeomorphism classes of tiled surface from \cite[Thm.~4.10]{PPP19}, by associating relations generating $I/J$ with the labels on the surface. 
\end{proof}

We refer the reader to \cite[Definition 4.6]{PPP19} for a description of the surface associated to a given locally gentle algebra. Here the authors use a pair of dual dissections of the surface. Our partial triangulation is their green dissection, and our marked points in $\sR$ are the vertices of their red dissection, the dual of $\sP$. 
    
\begin{figure}[ht!]
   \centering
    \includegraphics[scale=.55]{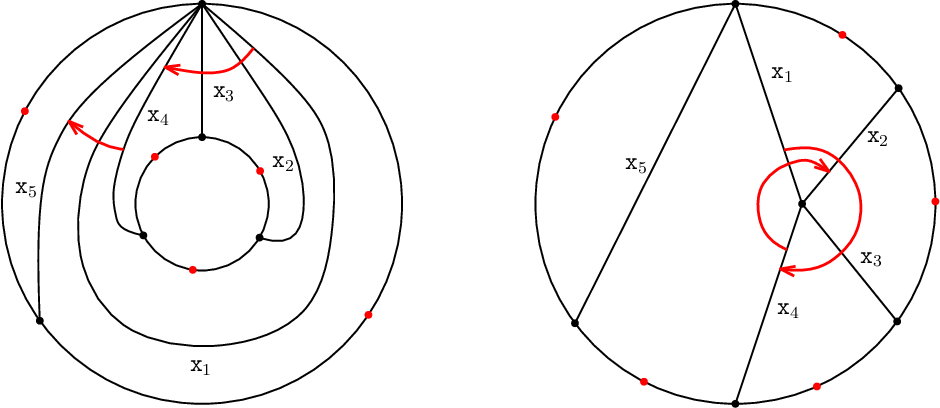}
    \caption{Non-homeomorphic surfaces associated to the same string algebra.}
    \label{fig:unique-surface}
\end{figure}

\begin{remark}\label{rem:choice-signs}
We highlight the following features of Theorem \ref{thm:model-algebra}.
\begin{enumerate}
    \item Given a string algebra $A$, the choice of an associated locally gentle algebra corresponds to the choice of sign functions $\sigma$ and $\epsilon$. Fix a choice, let $B$ be the corresponding locally gentle algebra, and assign to $B$ the same choice of sign functions. Then we have that $vw$ is a string in $B$ if and only if $\sigma (w) = -\epsilon (v)$. In particular, if $a,b$ are arrows of $A$ such that $ab$ corresponds to a label in the tiled surface associated to $B$, then $ab$ is not a relation in $B$ and so we have $\sigma (b) = - \epsilon (a)$. This observation will be used in the proofs of Theorem~\ref{thm:permissiblearcsstrings} and Theorem~\ref{thm:pivot-move-irreducible}.  
    \item The number of boundary components of the tiled surface associated to a locally gentle algebra is given in terms of perfect matchings of a graph whose vertices are so called blossom vertices, see~\cite[Remark 4.11 (vi)]{PPP19} for more details. It thus follows that a string algebra $A$ is a labelled tiling algebra of a surface with no boundary if and only if every vertex of the quiver of $A$ has valency four. 
    \item The existence of different choices of locally gentle algebras associated to a string algebra means that the surface associated to a string algebra is not unique in general. Example~\ref{ex:two-surfaces-one-algebra} illustrates this; see Figure~\ref{fig:unique-surface}. 
\end{enumerate}
\end{remark}

\begin{example}\label{ex:two-surfaces-one-algebra}
Consider the quiver $Q$ 
\[
\xymatrix{
4\ar[r]^{d}  & 1\ar[r]^{e}\ar[d]^{a}  & 5 \\
3\ar[u]^{c} & 2 \ar[l]^{b}
}
\] 
from Example~\ref{ex:string-gentle}
and the admissible ideal  
$I=\left<da, de, abc\right>$. The algebra $A=\kk Q/I$ is a string algebra, but not gentle. 

The two possible choices of locally gentle algebras are the finite-dimensional algebra $B_1 = \kk Q/J_1$ with $J_1 = \left<da\right>$ and the infinite dimensional algebra $B_2 = \kk Q/J_2$, with $J_2 = \left<de\right>$. The tiled surfaces corresponding to $B_1$ and $B_2$ are illustrated in Figure~\ref{fig:unique-surface} and are not homeomorphic to each other.

The algebra $A$ can then be associated to the surface of $B_1$ with labels corresponding to the relations $abc$ and $de$, and to the surface of $B_2$ with labels corresponding to the relations $abc$ and $da$. 
\end{example} 

\section{The module category of a labelled tiling algebra}\label{sec:geom-model}
In this section we will describe the module category of a string algebra via surfaces. 

Recall that $(\sS, \sM \cup \sR, \sP, \sL)$ denotes a labelled tiling, where the partial triangulation $\sP$ divides $\sS$ into a collection of polygonal tiles with at most one boundary segment in $(\sS, \sM)$ each. Each of these tiles has a marked point in $\sR$ which lies in the boundary segment if it exists, otherwise it is a puncture (see Figure~\ref{fig:tiles-red}). Recall also that $A_{\sP,\sL}$ denotes the corresponding labelled tiling algebra and any string algebra can be realised in this way.

\subsection{Permissible arcs and closed curves}

We will now adapt the definition of permissible arcs and closed curves from  \cite{BCS21} to any labelled tiling algebra. In particular, we will extend this definition to consider permissible arcs whose endpoints lie in $\sR$. These haven't been considered in~\cite{BCS21}, but they have been considered in~\cite{Chang} and as we will see in Subsection~\ref{sec:morphisms}, sometimes it is useful to consider these arcs instead. 

We first need to set up some definitions and notation.

\begin{definition}
    Let $\gamma, \gamma'$ be two curves in $(\sS,\sM \cup \sR)$.
    \begin{enumerate}
        \item We say {\it $\gamma$ consecutively crosses $\arc{a}_1, \arc{a}_2 \in \sP$} if $\gamma$ crosses $\arc{a}_1$ and $\arc{a}_2$ at the points $p_1$ and $p_2$, and the segment of $\gamma$ between the points $p_1$ and $p_2$ does not cross any other arc in $\sP$. 
        \item We denote by $I (\gamma, \gamma')$ the minimal number of transversal intersections of representatives of the isotopy classes of $\gamma$ and $\gamma'$ in the interior of $\sS \setminus (\sM \cup \sR)^\circ$. 
        \item The {\it intersection vector} $I_\sP (\gamma)$ of $\gamma$ with respect to $\sP$ is the vector $(I(\gamma, \arc{a}_i))_{\arc{a}_i \in \sP}$.
        \item The {\it intersection number $|I_\sP (\gamma)|$ of $\gamma$ with respect to $\sP$} is given by $\sum_{\arc{a} \in \sP} I(\gamma, \arc{a})$.
        \item $\gamma$ is said to be a {\it zero arc} if $I_\sP (\gamma)$ is the zero vector, i.e. $\gamma$ does not intersect any arc in $\sP$ except possibly at an endpoint. 
    \end{enumerate}
\end{definition}

Note that we might have $I(\gamma, \gamma') = 0$ for arcs $\gamma$, $\gamma'$ sharing an endpoint. In particular, if $\gamma \in \sP$, then $I(\gamma, \gamma') = 0$, for every arc $\gamma' \in \sP$. 

It will be convenient for us to realise each non-zero arc $\gamma$ as a concatenation of segments, defined by the consecutive crossings with $\sP$. As such, we will often denote it by $\gamma=\gamma_1 \cdots \gamma_r$, with $r \geq 2$. The segments $\gamma_1$ and $\gamma_r$ are called {\it end segments} and the remaining ones are called {\it interior segments}.  

\begin{definition}\label{def:M-R-arcs}
Let $\gamma$ be an arc in $(\sS, \sM \cup \sR)$.
\begin{enumerate}
\item $\gamma$ is an {\it $\sR$-arc} if its endpoints lie in $\sR$.
\item $\gamma$ is an {\it $\sM$-arc} if each endpoint of $\gamma$ lies in $\sM$ or it is a red point lying in a tile of type (i.b) (See Figure~\ref{fig:tiles-red}).
\end{enumerate}
\end{definition} 

We note that all arcs in \cite{BCS21} are what we now call $\sM$-arcs. Indeed, in the model used in~\cite{BCS21}, the only isolated marked points, i.e. the marked points not incident with the dissection $\sP$, were the ones lying in tiles of type (i.b).

We will see in Section~\ref{sec:morphisms} that $\sR$-arcs are useful to describe morphisms between string modules. On the other hand, $\sM$-arcs will be important in the classification of support $\tau$-tilting modules (Section~\ref{sec:support-tau-tilting}). For this classification, we make use of irreducible morphisms and the Auslander-Reiten translate. Therefore, we will consider $\sM$-arcs in their geometric descriptions in Sections~\ref{sec:pivot} and~\ref{sec:AR-translate}.

\begin{figure}[ht]
\includegraphics[height=3.3cm]{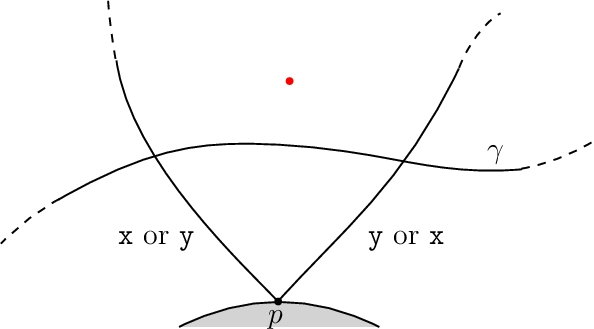}
\caption{A consecutive cross cuts an angle of $\sP$.}
\label{fig:permissible}
\end{figure}

\begin{definition}\label{def:permissible-arcs}
Let $\gamma$ be a curve in $(\sS, \sM \cup \sR)$.
\begin{compactenum}
\item $\gamma$ is called {\it permissible} if it satisfies the following conditions: 
\begin{compactenum}
\item If $\gamma$ consecutively crosses two (possibly not distinct) arcs $\arc{x}$ and $\arc{y}$ of $\sP$, then $\arc{x}$ and $\arc{y}$ have a common endpoint $p \in M$, and $\gamma$ cuts the corresponding angle at $p$; see Figure~\ref{fig:permissible}. 
\item No segment of $\gamma$ cuts a label in $\sL$. 
\end{compactenum} 
\item $\gamma$ is a {\it permissible $\sM$-arc} ({\it $\sR$-arc}, resp.) if it is both permissible and an $\sM$-arc ($\sR$-arc, resp.). 
\item If $\gamma$ is a closed curve, define $\gamma^n\colon [0,n] \to \sS$, where $n \in \mathbb{N}$, by $\gamma^n (x) \coloneqq \gamma(x-\lfloor{x}\rfloor)$ for $x \in [0,n]$. 
\item $\gamma$ is a {\it permissible closed curve} if it is a closed curve which satisfies (1)(a) and such that no segment of $\gamma^n$ cuts a label in $\sL$ for any $n \in \mathbb{N}$. 
\end{compactenum}
\end{definition}

We note that this definition generalises that in~\cite{BCS21}. The condition in~\cite{BCS21} regarding the winding number is not used here. In fact, this condition was used to describe a specific permissible $\sM$-arc used in the geometric description of irreducible morphisms and the AR-translate. However, we can give a different description of these $\sM$-arcs, without referring to the winding number (see Definition~\ref{def:pivot-move-t}).

\begin{figure}[ht!]
\includegraphics[scale=0.6]{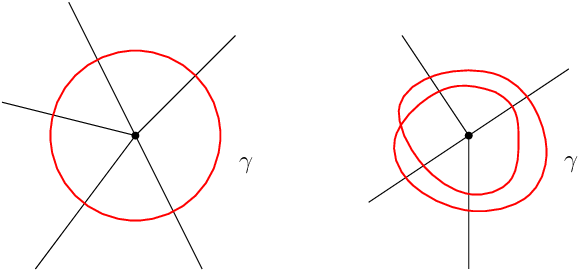}
\caption{Curves wrapping around a single puncture are not permissible.}
\label{fig:non-permissible-curve}
\end{figure}

\begin{remark}
    It follows from the definition that no closed curve which wraps around a puncture in $\sM$ (see Figure~\ref{fig:non-permissible-curve}) is permissible, because every puncture in $\sM$ has a label. 
\end{remark}

We will regard permissible arcs up to equivalence, which requires more than just isotopy. The following definition of equivalence of permissible arcs coincides with that in \cite{BCS21}. 

\begin{definition}\label{def:equiv-arcs}
Two permissible arcs $\gamma$ and $\gamma'$ in $(\sS, \sM \cup \sR)$ are said to be {\it equivalent} if we can write $\gamma = \gamma_1 \cdots \gamma_r$, $\gamma'= \gamma'_1 \cdots \gamma'_r$, where $\gamma_i$ and $\gamma'_i$ cut the same angle, for each $2 \leq i \leq r-1$.  
\end{definition}

\begin{figure}[ht!]
\includegraphics[height=4cm]{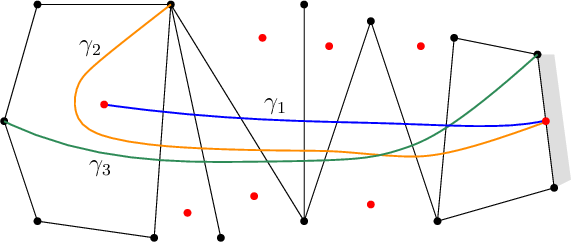}
\caption{These three arcs are permissible and equivalent to each other. The arc $\gamma_1$ (in blue) is the unique $\sR$-arc in this equivalence class. The arc $\gamma_3$ is an $\sM$-arc and $\gamma_2$ is neither an $\sM$-arc nor an $\sR$-arc.}
\label{fig:equivalence}
\end{figure}

Note that the end-segments do not play a role in the definition of equivalence of permissible arcs. The point of this definition is that equivalent permissible arcs correspond to the same string (see Theorem~\ref{thm:permissiblearcsstrings}). Moreover, each equivalence class of permissible arcs contains exactly one permissible $\sR$-arc up to isotopy (see Figure~\ref{fig:equivalence} for an illustration).

\begin{theorem}\label{thm:permissiblearcsstrings} 
Let $A$ be a string algebra and $(\sS, \sM \cup \sR, \sP, \sL)$ a labelled tiling associated to $A$. There are one-to-one correspondences between: 
\begin{compactenum}
\item equivalence classes of permissible arcs in $(\sS, \sM \cup \sR)$ and the strings of $A = A_{\sP, \sL}$.
\item isotopy classes of permissible closed curves $c$ in $(\sS, \sM \cup \sR)$ with $|I_\sP (c)| \geq 2$ and powers of bands of $A = A_{\sP, \sL}$.
\end{compactenum}
\end{theorem}
\begin{proof} 
To prove this we use a similar strategy as in the gentle case (see \cite[Theorems 3.8 and 3.9]{BCS21}). First, we associate the zero string to any zero arc. Note that this includes boundary segments and arcs in $\sP$. For non-zero permissible arcs with intersection number with respect to $\sP$ greater or equal to $2$, each consecutive crossing of a curve with $\sP$ corresponds to an arrow, and since all intersections with $\sP$ are considered to be transversal, we obtain a reduced walk. Now suppose that this reduced walk contains a direct path $\alpha_1 \cdots \alpha_r$ which is a relation. Again by transversality of the intersections with $\sP$ and by condition 1(a) in Definition~\ref{def:permissible-arcs}, this relation is not a relation in $A_{\sP}$. Therefore, the path $\alpha_1 \cdots \alpha_r$ corresponds to a label, contradicting condition 1(b) in Definition~\ref{def:permissible-arcs}. 

The correspondence between (permissible) arcs with only one intersection with $\sP$ and the trivial strings of $A_{\sP,\sL}$ is as follows. Let $\sigma$ and $\epsilon$ be the sign functions used to define $A_\sP$ (see proof of Theorem~\ref{thm:model-algebra}). 

Given an arc $\gamma = \gamma_1 \gamma_2$ with only one intersection with $\sP$, we associate to it the trivial string $1_x^\pi$, where $x$ and $\pi$ are as in Figure~\ref{fig:trivialstringfromarc}. Note that not all of the four arrows $a_1, b_1, a_2, b_2$ have to exist, but if $A$ is not an algebra of type $A_1$, at least one of them does exist. 

\begin{figure}[ht!]
\includegraphics[scale=0.6]{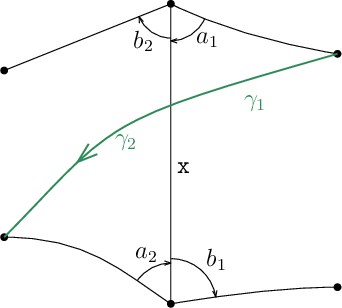}
\caption{The arc $\gamma = \gamma_1 \gamma_2$ corresponds to the trivial string $1_x^\pi$, where $\pi = \epsilon (a_1) = \sigma (b_1) = - \epsilon (a_2) = - \sigma (b_2)$.}
\label{fig:trivialstringfromarc}
\end{figure}

Conversely, given a string $1_x^\pi$, we associate to it an arc which intersects $\sP$ once at arc $\arc{x}$, and we orient this arc in such a way that $\gamma$ satisfies the equalities described in Figure~\ref{fig:trivialstringfromarc}. Note that these equalities hold due to Remark~\ref{rem:choice-signs} (1).

Finally, given a permissible closed curve $\gamma$, pick a crossing $c$ with $\sP$, let $c= \gamma (0)$ and pick an orientation of $\gamma$. The corresponding power of a band $b^k$, for some $k \geq 1$, is given by the arrows corresponding to the angles that $\gamma$ cuts from $\gamma (0)=c$ to $\gamma (1)=c$. 
\end{proof}

\begin{remark}
The correspondence between permissible arcs and strings in Theorem~\ref{thm:permissiblearcsstrings} is such that:
\begin{itemize}
    \item the inverse $\gamma^{-1}$ of a permissible arc $\gamma$ corresponds to $w_\gamma^{-1}$,
    \item the intersection vector of a permissible arc $\gamma$ gives the dimension vector of the corresponding string module. 
\end{itemize}
\end{remark}

Given a permissible curve $\gamma$ in $(\sS, \sM \cup \sR)$, we will denote by $w_\gamma$ the corresponding string of $A$. If the curve $\gamma$ is an arc, the corresponding string module will be denoted by $M_\gamma$.

\begin{remark}\label{rem:permissible}
It follows from Theorem~\ref{thm:permissiblearcsstrings} that a labelled tiling algebra is of finite representation type if and only if the intersection number of any permissible closed curve with respect to $\sP$ is at most one. We note that in the gentle case, it is enough to check this property for simple closed curves, cf.~\cite[Corollary 3.10]{BCS21}. However, this is not the case for string algebras in general. For example, the string algebra given in Example~\ref{ex:labeled} is of infinite representation type, as it has two bands: $b_1 \coloneqq a_7a_4a_6a_8^{-1}$ and $b_2 \coloneqq a_2 a_5^{-1} a_6 a_8^{-1} a_7 a_3^{-1}$. These bands correspond to the closed curves in Figure~\ref{fig:bands}. Note that in this example there are no simple closed curves which are permissible. 

\begin{figure}
\centering
\includegraphics[scale=.6]{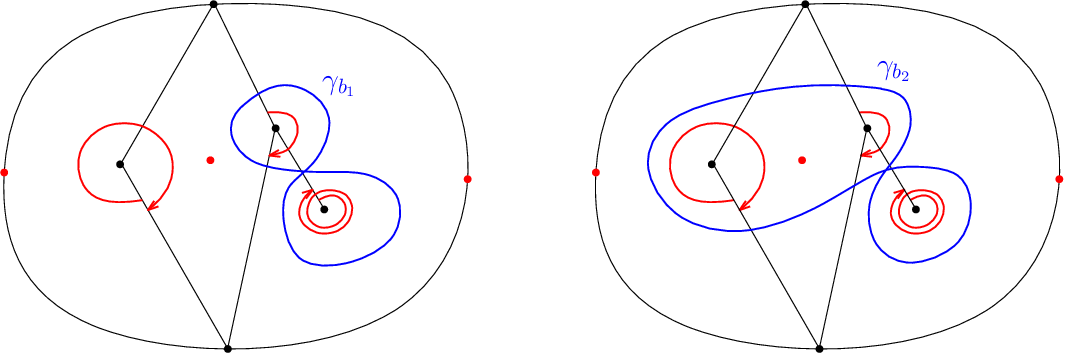}
\caption{Permissible closed curves corresponding to $b_1$ (on the left) and to $b_2$ (on the right).}
    \label{fig:bands}
\end{figure}
\end{remark}

\subsection{Morphisms via crossings and intersections}\label{sec:morphisms}

In this section, we provide a combinatorial description of morphisms between string modules. We will make use of permissible $\sR$-arcs, as they permit us to give a simple description of morphisms via crossings and intersections. Such a description has been given in~\cite[Section 2.5]{Chang} in the case when the algebra is gentle. The proof given in~\cite{Chang} is via an embedding of the module category of a gentle algebra in the bounded derived category and therefore cannot be generalised to string algebras. Our proof is different and it relies on string combinatorics, namely Theorem~\ref{thm:homs}. 

\begin{figure}[ht!]
 \centering
\includegraphics[scale=.8]{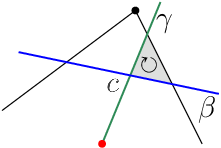}
\hskip.5cm
\includegraphics[scale=.8]{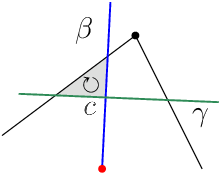}
\hskip.5cm
\includegraphics[scale=.8]{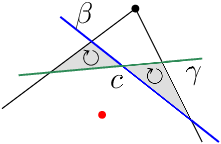}
\hskip.5cm
\includegraphics[scale=.9]{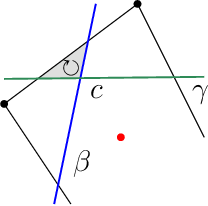}
\caption{A crossing between $\sR$-arcs, from $\beta$ to $\gamma$.}
\label{fig:crossfrombetatogamma}
\end{figure}

\begin{definition}\label{def:crossing-R}
    Let $\gamma$ and $\beta$ be two permissible $\sR$-arcs, and suppose there is a crossing $c$ between $\gamma$ and $\beta$ in the interior of $\sS$. Then there must be an arc $\arc{x}$ in $\sP$ crossed by both arcs $\gamma$ and $\beta$. The arcs $\arc{x}, \gamma$ and $\beta$ enclose a simply connected triangle $\delta$ with $c$ as one of its vertices. We say the crossing $c$ is a {\it crossing from $\beta$ to $\gamma$} if $\gamma$ follows $\beta$ in the clockwise order around $\delta$. See Figure~\ref{fig:crossfrombetatogamma}. 
\end{definition}

Note that a given crossing is incident with at most two triangles which are opposite to each other. Thus both triangles give rise to the same orientation of the crossing. Moreover, the orientation is independent of the choice of representatives of the arcs in their isotopy classes. Therefore, the orientation of the crossing is well defined.

\begin{figure}[ht!]
 \centering
\includegraphics[scale=.7]{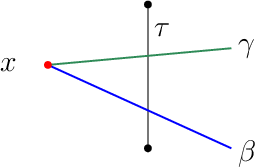}
\hskip.6cm
\includegraphics[scale=.8]{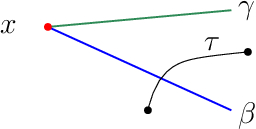}
\hskip.6cm
\includegraphics[scale=.7]{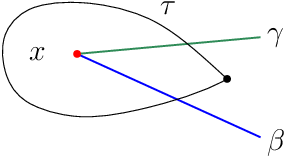}
\caption{Left: Intersection from $\beta$ to $\gamma$ traversed by $\sP$. Middle: Intersection from $\beta$ to $\gamma$ not traversed by $\sP$. Right: Intersection from $\beta$ to $\gamma$ not traversed by $\sP$ and intersection from $\gamma$ to $\beta$ traversed by $\sP$.}
\label{fig:intersection-traversed}
\end{figure}

\begin{definition}\label{def:crossings-R}
  Let $\gamma$ and $\beta$ be two permissible $\sR$-arcs.
  \begin{enumerate}
      \item We say there is an {\it intersection from $\beta$ to $\gamma$} if there are end-segments $\beta_0$ and $\gamma_0$ of $\beta$ and $\gamma$ respectively sharing an endpoint $x \in \sR$ and such that $\gamma_0$ follows $\beta_0$ in the anticlockwise direction around $x$.
      \item A {\it trivial intersection from $\beta$ to $\beta$} is an intersection from $\beta$ to $\beta$ in which the two end-segments coincide. 
      \item An intersection from $\beta$ to $\gamma$ is {\it traversed by $\sP$} if there is an arc $\tau$ in $\sP$ which crosses both $\beta$ and $\gamma$ and such that $\beta$, $\gamma$ and $\tau$ enclose a simply connected triangle. See Figure~\ref{fig:intersection-traversed} for an illustration. 
  \end{enumerate}
\end{definition}

We note that our definition of intersections traversed by $\sP$ is motivated by \cite{BCSGS24}, where these are described in terms of certain angles at red points. 

Denote by $\sC(\beta,\gamma)$ the number of crossings from $\beta$ to $\gamma$ and by $\sI_{\sR,\sP}(\beta,\gamma)$ the number of intersections from $\beta$ to $\gamma$ traversed by $\sP$.

\begin{theorem}\label{thm:R-arcs-morphisms}
    Let $\beta$ and $\gamma$ be permissible $\sR$-arcs, and $M_\beta, M_\gamma$ the corresponding string modules. Then $\dim \Hom(M_\beta, M_\gamma) = \sC(\beta,\gamma)+\sI_{\sR,\sP} (\beta,\gamma)-\delta_{\beta,\gamma}$.
\end{theorem}
\begin{proof}
The dimension of $\Hom(M_\beta, M_\gamma)$ is given by the number of overlaps in $\Fac (\beta) \cap \Sub (\gamma)$.
Let $w_e$ be such an overlap, and assume first that $M_\beta \not\cong M_\gamma$.

{\bf Case 1:} $w_e$ is not a prefix or suffix of one of the strings, say $w_\beta$. Then $w_\beta = w_\ell a^{-1} w_e b w_r$, for some $a, b \in Q_1$. On the other hand we have $w_\gamma = v_\ell c w_e d^{-1} v_r$, where we may have $v_\ell c=0$ and/or $d^{-1} v_r = 0$. 
In any case, we must have a crossing from $\beta$ to $\gamma$, as Figure~\ref{fig:e-not-prefix-suffix} illustrates. 

{\bf Case 2:} $w_e$ is a prefix or suffix of both strings $w_\beta$ and $w_\gamma$. 

{\bf Case 2.1:} $w_e$ is non-trivial. We can assume, without loss of generality, that either $w_e$ is a prefix of both $w_\beta$ and $w_\gamma$ or $w_e$ is a prefix of $w_\beta$ and suffix of $w_\gamma$. In the first case, there is an intersection from $\beta$ to $\gamma$ traversed by $\sP$. In the second case, we can assume $w_\beta = w_e b w_r$ and $w_\gamma = v_\ell c w_e$, with $b, c \in Q_1$, as otherwise we would be in the first case situation. As Figure~\ref{fig:e-suffix-and-prefix} illustrates, we must have a crossing from $\beta$ to $\gamma$. 

\begin{figure}[ht!]
 \centering
\includegraphics[scale=.8]{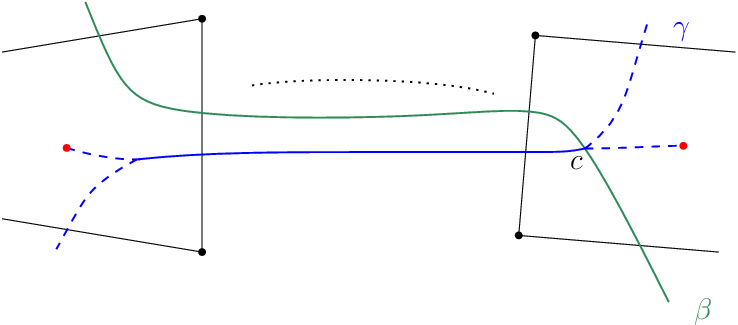}
\caption{Case where $w_e$ is not a suffix or prefix of one of the strings, say $w_\beta$.}
\label{fig:e-not-prefix-suffix}
\end{figure}

\begin{figure}[ht!]
\centering
\includegraphics[scale=.8]{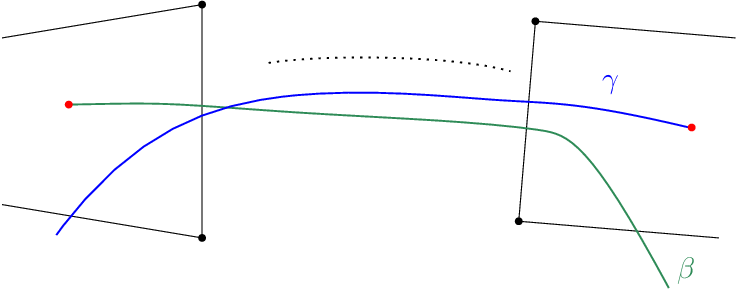}
\caption{Case where $w_e$ is a prefix of $w_\beta$ and a suffix of $w_\gamma$.}
\label{fig:e-suffix-and-prefix}
\end{figure}

{\bf Case 2.2:} $w_e$ is trivial, i.e. $w_e=1_x^{\pi}$, with $x\in Q_0$ and $\pi \in \{+,-\}$. Without loss of generality, we can write $w_\beta = b w_r$ and $w_\gamma = d^{-1} v_r$, with $b \in Q_1$ (resp. $d^{-1} \in Q_1^{-1}$) if $w_\beta$ (resp. $w_\gamma$) is non-trivial, such that $s(b) = t(d) = x$.  Since $M_\beta \not\cong M_\gamma$, at least one of $w_\beta$ and $w_\gamma$ is non-trivial. Without loss of generality, assume the latter. 

If $w_\beta$ is trivial or $db=0$ is not a label, then we have an intersection from $\beta$ to $\gamma$ traversed by $\sP$ (see Figure~\ref{fig:DB-label-or-B-empty}). Otherwise, either $db$ is a label or $db \neq 0$, in which case we have a crossing from $\beta$ to $\gamma$, as Figure~\ref{fig:DB-label-or-DB-nonzero} illustrates. 

\begin{figure}[ht!]
 \centering
\includegraphics[scale=.8]{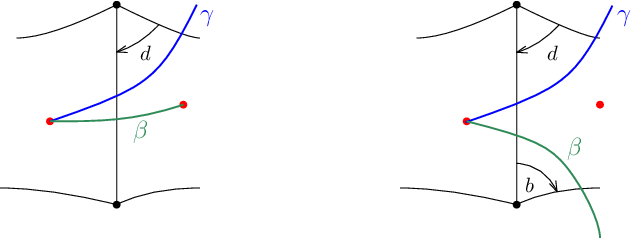}
\caption{Case where $w_\beta$ is trivial or $db = 0$ but not a label.}
\label{fig:DB-label-or-B-empty}
\end{figure}

\begin{figure}[ht!]
 \centering
\includegraphics[scale=.8]{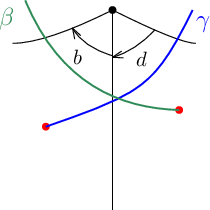}
\caption{Case where $db$ is a label or $db \neq 0$.} 
\label{fig:DB-label-or-DB-nonzero}
\end{figure}

If $M_\beta \cong M_\gamma$, then there are two intersections from $\beta$ to itself, at both end-segments. Both intersections are trivial and traversed by $\sP$. We can associate these trivial intersections to the overlap $w_e=w_\beta$. The corresponding morphism is the identity morphism. 

In summary, we associated one crossing or one non-trivial intersection from $\beta$ to $\gamma$ traversed by $\sP$ to a morphism from $M_\beta$ to $M_\gamma$ which is not the identity morphism, and we associated the two trivial intersections from $\beta$ to $\gamma$ to the identity morphism, in the case when $M_\beta \cong M_\gamma$. 

Conversely, let $c$ be a crossing from $\beta$ to $\gamma$, and let $\alpha_i$ and $\beta_j$ be the segments of $\alpha$ and $\beta$, respectively, where the crossing takes place. 
We can then consider the maximal sequences of consecutive segments:
\[
(\gamma_{i-s}, \ldots, \gamma_i, \gamma_{i+1}, \ldots, \gamma_{i+r}) \text { and } (\beta_{j-s}, \ldots, \beta_j, \beta_{j+1}, \ldots, \beta_{j+r})
\]
of $\gamma$ and $\beta$, respectively, which read off the same string (see Figure~\ref{fig:e-corresponding-to-crossing} for an illustration). Denote such a string by $w_e$. Note that $w_e$ can be a trivial string, but it is always non-zero, since there is at least one arc in $\sP$ that is crossed both by $\beta$ and $\gamma$. Given the crossing $c$ is from $\beta$ to $\gamma$, we must have $w_\beta = w_\ell a^{-1} w_e b w_r$ and $w_\gamma = v_\ell c w_e d^{-1} v_r$, where at least one of $w_\ell a^{-1}, v_\ell c$ is non-zero and at least one of $b w_r, d^{-1} v_r$ is non-zero. It then follows from Theorem~\ref{thm:homs}, that there is a morphism from $M_\beta$ to $M_\gamma$ corresponding to the overlap $w_e$. 

\begin{figure}[ht!]
 \centering
 \includegraphics[scale=.85]{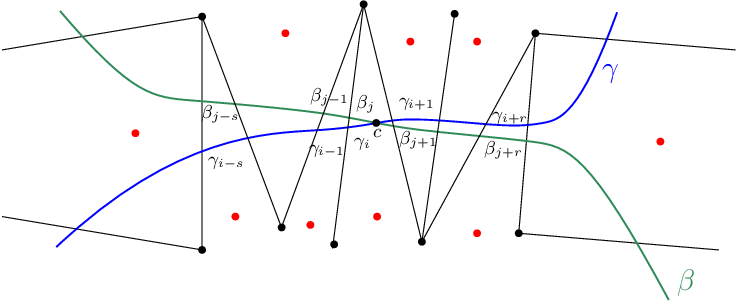}
\caption{The maximal sequence of segments reading off the same string.}
\label{fig:e-corresponding-to-crossing}
\end{figure}

Now, let $c$ be an intersection from $\beta$ to $\gamma$ traversed by an arc in $\sP$. Suppose this intersection is non-trivial.

{\bf Claim:} There are sequences of consecutive segments 
\[
(\beta_1, \ldots, \beta_r, \beta_{r+1}) \text{ and } (\gamma_1, \ldots, \gamma_r, \gamma_{r+1})
\]
of $\beta$ and $\gamma$ respectively, for some $r\ge 1$, such that $s(\beta_1) = s(\gamma_1)=c$,  $\beta_1 \ldots \beta_r$ and $\gamma_1 \ldots \gamma_r$ read off the same string $w_e$, and $\beta_{r+1}, \gamma_{r+1}$ read off different strings (one of which may be a zero string). 

{\it Proof of the claim.} If $M_\beta \not\cong M_\gamma$, this is clear. So suppose $M_\beta \cong M_\gamma$. Then a non-trivial intersection from $\beta$ to $\gamma=\beta$ 
would imply that both end-segments of $\beta$ meet at a point $x \in \sR$. Suppose for a contradiction that $w_\beta = w_\beta^{-1}$. Then $w_\beta = a_1 \cdots a_n$ is non-trivial, with either $a_{\lceil \frac{n}{2} \rceil} = a_{\lceil \frac{n}{2} \rceil}^{-1}$, if $n$ is odd, or $a_{\frac{n}{2}} = a_{\frac{n}{2}+1}^{-1}$, if $n$ is even, both a contradiction. Therefore $w_\beta \neq w_\beta^{-1}$, and the claim follows. $\blacksquare$

Since the intersection is from $\beta$ to $\gamma$ it follows that $w_e$ is a factor substring of $w_\beta$ and an image substring of $w_\gamma$. Hence, by Theorem~\ref{thm:homs}, we have a morphism from $M_\beta$ to $M_\gamma$ corresponding to the overlap $w_e$. 

Finally, if the intersection $c$ is trivial, then it corresponds to the identity morphism. But for each arc $\beta$ there are two such trivial intersections, so we are counting the identity morphism twice. 
\end{proof}

\begin{remark}
Let $A$ be a string module, $A= B/I$, where $B$ is a locally gentle algebra, and let $(\sS, \sM \cup \sR, \sP, \sL)$ be the corresponding labelled tiling. Given two string $A$-modules $M_\beta$ and $M_\gamma$, with $\sR$-arcs $\beta$ and $\gamma$ then $\sC(\beta, \gamma)$ is the number of overlap extensions from $M_\beta$ to $M_\gamma$ in $\mod B$ (see \cite{CPS21, Chang} for more details). 
\end{remark}

\begin{example}
   Consider the string algebra $A$ given by the labelled tiling in Figure~\ref{fig:label-angle-red} 
   (see Figure~\ref{fig:quiver-tiling-red} for the quiver of this tiling). 
 
Let $w=a_7^{-1}a_8$ and $v=a_8 a_6^{-1}$. Figure~\ref{fig:w-v-crossings} illustrates the five crossings and intersections, $c_1, \ldots, c_5$, between the corresponding permissible $\sR$-arcs. 

\begin{figure}[ht!]
 \centering
 \includegraphics[scale=1]{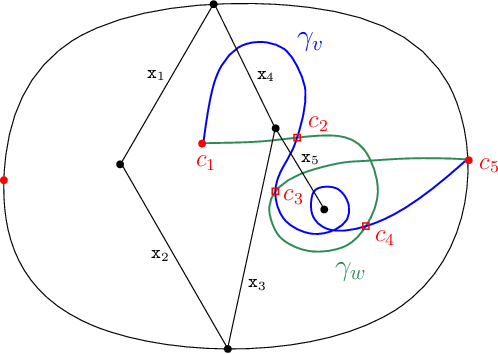}
\caption{The crossings and intersections between the $\sR$-arcs $\gamma_w, \gamma_v$ corresponding to strings $w$ and $v$.}
\label{fig:w-v-crossings}
\end{figure}

The intersections $c_1$ and $c_5$ are not traversed by $\sP$. The crossing $c_2$ is a crossing from $\gamma_w$ to $\gamma_v$ corresponding to the admissible pair $(a_7^{-1}1_5^\pi a_8,a_81_5^\pi a_6^{-1})$ in $\Fac (w) \times \Sub (v)$, with overlap $1_5^\pi$, where $\pi = -\sigma(a_8)$.  
The crossing $c_3$ is from $\gamma_w$ to $\gamma_v$ and it corresponds to the admissible pair 
$(a_7^{-1}a_8,a_8a_6^{-1})$ in $\Fac (w) \times \Sub (v)$ with overlap 
$a_8$.  
Finally the crossing $c_4$ is a crossing from $\gamma_v$ to $\gamma_w$, and it corresponds to the admissible pair 
$(1_5^\pi a_8a_6^{-1},a_7^{-1}a_8 1_5^\pi)$ 
in $\Fac (v^{-1}) \times \Sub (w)$, with overlap $1_5^\pi$, where $\pi = -\sigma (a_8)$. 

The dimension of $\Hom (M(w),M(v)) = 2$, the dimension of $\Hom (M(v),M(w)) = 1$, $I_\sP (v,w) = 0 = I_\sP (w,v)$, $C (w,v) =2$ and $C(v,w) =1$. 
\end{example}
 
\subsection{Pivot elementary moves}\label{sec:pivot}

In order to describe the Auslander-Reiten translate of a string module combinatorially, we first need to give a combinatorial description of the irreducible morphisms starting at string modules. We will give this description in terms of $\sM$-arcs instead of $\sR$-arcs because it allows us to give a more explicit classification of support tau-tilting modules (see Section~\ref{sec:support-tau-tilting}). 

The description of irreducible morphisms is split into two steps. The first step consists of choosing a specific representative in the equivalence class of permissible arcs.

\begin{figure}[ht!]
    \centering
    \includegraphics[scale=.6]{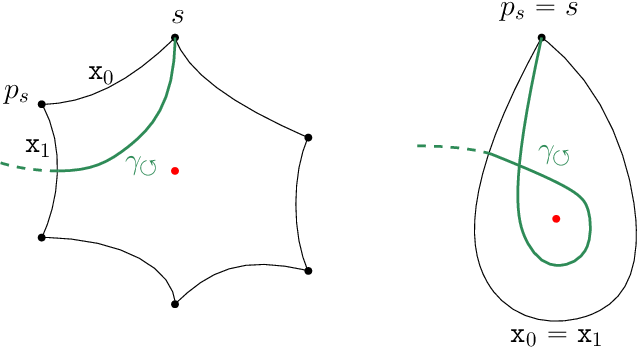}
    \caption{The end segments of the representative $\gamma_\circlearrowleft$.}
    \label{fig:representative}
\end{figure}

\begin{definition}
    Let $\gamma= \gamma_0 \gamma_1 \cdots \gamma_{t}$ be a permissible $\sM$-arc (or $\sR$-arc). Denote by  $\Delta_s$ (resp. $\Delta_t$) the tile in which $\gamma_0$ (resp. $\gamma_t)$ lies in. 
    
    The {\it anticlockwise-most $\sM$-arc representative} $\gamma_\circlearrowleft$ in the equivalence class of $\gamma$ is described as follows. The segments of $\gamma_\circlearrowleft$ which are interior segments are uniquely determined by the string defined by $\gamma$. The start segment (resp. terminal segment) of $\gamma_\circlearrowleft$ is such that the concatenation of $\gamma_\circlearrowleft$ with the side of $\Delta_s$ (resp. $\Delta_t$) incident with the starting point (resp. ending point) of $\gamma_\circlearrowleft$ in the anticlockwise direction gives rise to an arc which is not equivalent to $\gamma$. Figure~\ref{fig:representative} illustrates the end segments of $\gamma_\circlearrowleft$. 
\end{definition}

Note that we can dually define the {\it clockwise-most $\sM$-arc representative}, which will be denoted by $\gamma_\circlearrowright$.

We are now ready to describe the second step, which produces a new arc $f_t(\gamma)$ (resp. $f_s(\gamma)$) from $\gamma_\circlearrowleft$ by performing a pivot elementary move which fixes the endpoint (resp. starting point) of $\gamma_\circlearrowleft$.

\begin{figure}[t]
    \centering
    \includegraphics[scale=.65]{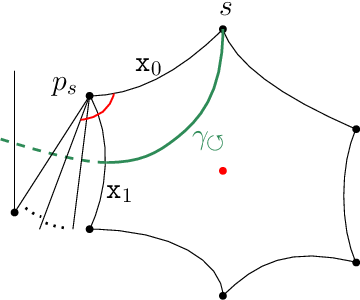}
    \hskip .2cm
    \includegraphics[scale=.65]{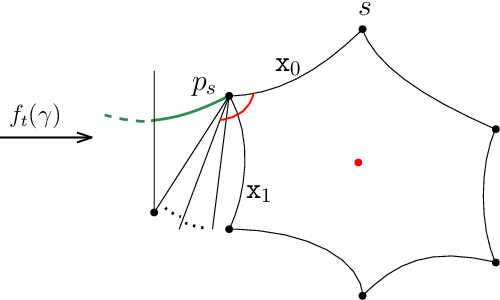}
    \caption{Pivot elementary move, case 1.}
    \label{fig:pivot-1}
\end{figure}

\begin{figure}[t]
    \centering
    \includegraphics[scale=.65]{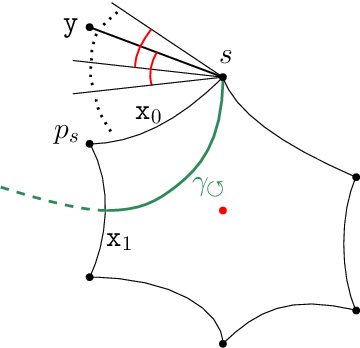}
    \hskip .2cm
    \includegraphics[scale=.65]{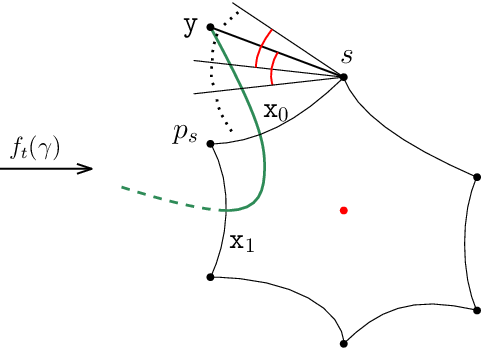}
    \caption{Pivot elementary move, case 2.}
    \label{fig:pivot-2}
\end{figure}

\begin{figure}[t]
    \centering
    \includegraphics[scale=.65]{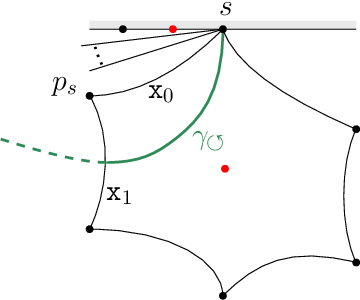}
    \hskip .2cm
    \includegraphics[scale=.65]{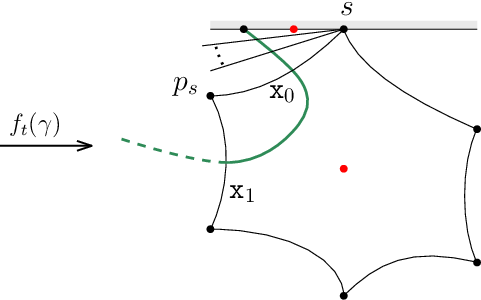}
    \caption{Pivot elementary move, case 3.}
    \label{fig:pivot-3}
\end{figure}

\begin{definition}\label{def:pivot-move-t} 
Let $\arc{x}_1$ be the first arc in $\sP$ crossed by $\gamma_\circlearrowleft$ and $\arc{x}_0$ be the side of $\Delta_s$ which precedes $\arc{x}_1$ with respect to the clockwise direction around the tile. Denote by $p_s$ the common endpoint of $\arc{x}_0$ and $\arc{x}_1$, and by $s$ the starting point of $\gamma_\circlearrowleft$. Denote by $f_t(\gamma)$ be the arc obtained from $\gamma_\circlearrowleft$ as follows:

{\bf Case 1.} There is a label at $p_s$ starting at the angle in $\Delta_s$ and such that it is {\it shorter than $\gamma_\circlearrowleft$}, i.e. it crosses fewer arcs in $\sP \setminus\{\arc{x}_0\}$ adjacent to $p_s$ than $\gamma_\circlearrowleft$ does. 

In this case, we define $f_t(\gamma)$ to be the concatenation of $\gamma_\circlearrowleft$ with $\arc{x}_0$ (see Figure~\ref{fig:pivot-1}). 

{\bf Case 2.} There is no label as in Case 1, but there are labels at $s$ starting after the angle in $\Delta_s$ in the clockwise direction around $s$. 

Consider the first label at $s$, when going around $s$ in the clockwise direction from $\arc{x}_0$, and let $\arc{y}$ be the last arc crossed by that label. We define $f_t (\gamma)$ to be the concatenation of $\gamma_\circlearrowleft$ with $\arc{y}$ (see Figure~\ref{fig:pivot-2}). 

{\bf Case 3.} There are no labels as in Cases 1 and 2. 

Note that this implies that $s$ is a marked point in the boundary of $\sS$, as otherwise, there would be a label at $s$, and we would be in the situation of Case 2. 

As such, we define $f_t(\gamma)$ to be the concatenation of $\gamma_\circlearrowleft$ with the boundary segment connecting $s$ with its predecessor in $\sM$ in the anticlockwise direction around the boundary component (see Figure~\ref{fig:pivot-3}). 

The arc $f_t(\gamma)$ is called the {\it pivot elementary move of $\gamma$ at its ending point}. The {\it pivot elementary move $f_s(\gamma)$ of $\gamma$ at its starting point} is defined dually.
\end{definition}

Note that if the algebra is gentle then we are always in Case 3, and we recover the definition in \cite{BCS21}. 

\begin{theorem}\label{thm:pivot-move-irreducible}
Let $w$ be a string, $\gamma$ be a permissible $\sM$-arc corresponding to $w$ and $M_\gamma$ be the corresponding string module. Each irreducible morphism starting at $M_\gamma$ is obtained by a pivot elementary move at an endpoint of $\gamma_\circlearrowleft$. 
\end{theorem}
\begin{proof}
Given a string $w$, recall that the orientation of any permissible arc associated to $w$ is well-defined if $w$ is non-trivial. In the case when $w$ is trivial, the orientation is determined using the fixed choice of sign functions; see Figure~\ref{fig:trivialstringfromarc}.

We will show that $f_t(\gamma_\circlearrowleft)$ is a permissible $\sM$-arc whose corresponding string is $f_t(w)$. The proof that $f_s(\gamma_\circlearrowleft)$ corresponds to $f_s(w)$ is similar. The result then follows from Theorem~\ref{thm:irreducible}.

{\bf Case a.} There is an arrow $a \in Q_1$ such that $aw$ is a string.

Then $f_t(w)$ is obtained by adding the hook of $a$ on the left of $w$, i.e. $p^{-1}_a a w$. Due to the orientation of $\gamma$, we have that $s(a) = x_0$, where $x_0$ is the vertex corresponding to the arc $\arc{x}_0$ in Definition~\ref{def:pivot-move-t}. Indeed this is easy to see when $w$ is non-trivial, and when $w$ is trivial, it must be of the form $1_{x_1}^{\epsilon (a)}$ since $aw$ is a string and so the orientation of $\gamma_\circlearrowleft$ is as in Figure~\ref{fig:trivialstring-case1}.

\begin{figure}[ht!]
    \centering
\includegraphics[scale=0.7]{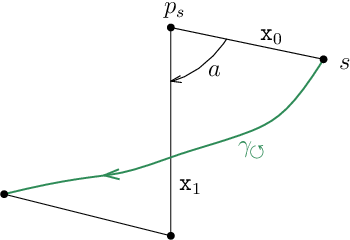}\caption{The arc $\gamma_\circlearrowleft$ corresponding to $1_{x_1}^{\epsilon (a)}$.}
    \label{fig:trivialstring-case1}
\end{figure}

Moreover, there cannot be a label at $p_s$ as in Case 1 of Definition~\ref{def:pivot-move-t}, as otherwise $a w$ would not be a string. We are therefore in either Case 2 or Case 3 of Definition~\ref{def:pivot-move-t}. In either case, the string corresponding to $f_t(\gamma)$ does indeed correspond to $p_a^{-1} a w$. In Case 2, the path $p_a$ is right maximal in the sense that there is $b \in Q_1$ with $s(b) = t(p_a)$ but $p_a b$ is a relation corresponding to a label. In Case 3, the path $p_a$ is right maximal because there is no $b \in Q_1$ with $s(b) = t(p_a)$. 

{\bf Case b.} There is no arrow $a$ such that $a w$ is a string.

This means that either $\arc{x}_0$ is a boundary segment in $(\sS, \sM)$, and so we are in Case 3 of Definition~\ref{def:pivot-move-t} or there is a label at $p_s$  as in Case 1 of Definition~\ref{def:pivot-move-t}. In both cases, we have that $f_t (\gamma)$ is the concatenation of $\gamma_\circlearrowleft$ with $\arc{x}_0$. 

This move removes the direct path $p$ at the start of $w$ which corresponds to a fan at $p_s$. Moreover, $f_t (\gamma)$ is either an arc with no intersections with $\sP$ if $w$ is direct or it is an arc whose corresponding string $w_r$ is such that $w = q_b b^{-1} w_r$, where $b \in Q_1$ and $q_b b^{-1}$ is the cohook of $b$. In either case, the string corresponding to $f_t(\gamma)$ is indeed $f_t(w)$, and we are done. 
\end{proof}

%
\subsection{Geometric interpretation of the Auslander-Reiten translate}\label{sec:AR-translate}

We are now ready to describe the (inverse of the) Auslander-Reiten translate of a string module in the geometric model. 

\begin{definition}
    Let $\gamma$ be a permissible arc, and consider the corresponding anticlockwise most $\sM$-arc representative $\gamma_\circlearrowleft$. Let $\arc{a}_s$ and $\arc{a}_t$ be the arcs of $\sP$ or boundary segments of $(\sS, \sM)$ such that $f_t(\gamma)$ (resp. $f_s(\gamma)$) is the concatenation of $\gamma_\circlearrowleft$ with $\arc{a}_s$ (resp. $\arc{a}_t$). Define $\tau^{-1} (\gamma)$ to be the concatenation of $\gamma_\circlearrowleft$ with the arcs $\arc{a}_s$ and $\arc{a}_t$. The arc $\tau (\gamma)$ is defined dually using $\gamma_\circlearrowright$. 
\end{definition}

The action of $\tau^{-1}$ (or $\tau$) on an arc can be seen as a generalisation of the ``classical'' description, in the Jacobian algebras case, of the Auslander-Reiten translate in terms of rotation of the arc. The following proposition generalises Theorem 3.9 and Corollary 3.17 of \cite{BCS21}. 

\begin{proposition}
Let $M$ be an indecomposable module over a labelled tiling algebra $A$ associated to $(\sS, M, \sP, \sL)$ and $\gamma$ a curve corresponding to $M$.
\begin{compactenum}
\item $M$ is injective at vertex $v$ if and only if $\tau^{-1} (\gamma) = \arc{v} \in \sP$. 
\item $M$ is a non-injective string module if and only if $\tau^{-1} (\gamma)$ has non-zero intersection with $\sP$, in which case the module corresponding to $\tau^{-1} (\gamma)$ is $\tau^{-1}(M)$. 
\item If $M$ is a non-injective string module and $f_s(\gamma)\neq 0$ (resp. $f_t(\gamma) \neq 0$), then $\tau^{-1} (\gamma) = f_t(f_s(\gamma))$ (resp. $\tau^{-1} (\gamma) = f_s(f_t(\gamma))$.
\item If $M$ is a band module then $\gamma$ corresponds also to $\tau^{-1} (M)$. 
\end{compactenum}
\end{proposition}
\begin{proof}
Suppose $M$ is a string module. Then $M$ is injective if and only if the arc $\gamma_\circlearrowleft$ meets only two fans which are arranged as in Figure~\ref{fig:injective}. Note that $\arc{a}_s$ (resp. $\arc{a}_t$) is either the last arc crossed by the first label at $p_s$ (resp. $p_t$) which starts at or after $\arc{v}$ in the anticlockwise direction around $p_s$  (resp. $p_t$) or $\arc{a}_s$ (resp. $\arc{a}_t$) is a boundary segment in $(\sS, \sM)$ if there is no such label. In any case, condition $(1)$ follows.

\begin{figure}[ht!]
    \centering
    \includegraphics[scale=0.7]{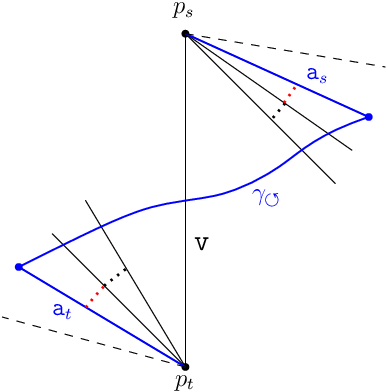}
    \caption{The arc $\gamma_\circlearrowleft$ corresponding to $I_v$.}
    \label{fig:injective}
\end{figure}

In particular, if $M$ is injective, then $\tau^{-1} (\gamma)$ does not intersect $\sP$ traversely. In order to prove condition $(2)$ it remains to show the forward implication. So suppose $M$ is a non-injective string module. Lemma~\ref{lem:non-injective} and Theorem~\ref{thm:pivot-move-irreducible} imply that $f_s(\gamma)\neq 0$ or $f_t(\gamma)\neq 0$. Assume, without loss of generality, that $f_s(\gamma)\neq 0$. By Theorem~\ref{thm:irreducible}, $\tau^{-1} (M)$ corresponds to the arc $f_t(f_s(\gamma))$. So it is enough to show that $\tau^{-1} (\gamma) = f_t(f_s(\gamma))$, which would also prove condition $(3)$. We do this now.

{\bf Case a.} $f_s(w_\gamma)$ is obtained by adding a hook $a^{-1} p_a$ on the right.

In this case, we can either add a hook $p_c^{-1} c$ on the left of $f_s(w_\gamma)$ or remove a cohook $q_c c^{-1}$ on the left of $f_s (w_\gamma)$. The latter is always possible since $f_s(w_\gamma) = w a^{-1} p_a$ is not a direct string. Figure~\ref{fig:add-hook-on-right} illustrates the different possibilities. In any of the cases, we have $\tau^{-1} (\gamma) = f_t(f_s(\gamma))$. 

\begin{figure}[ht!]
    \centering
\includegraphics[scale=0.6]{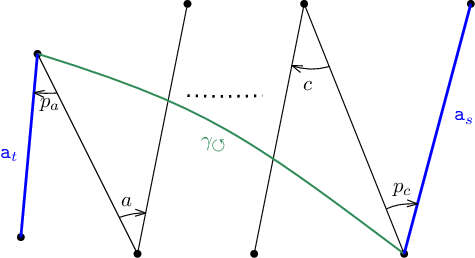}
    \hskip.5cm
\includegraphics[scale=0.6]{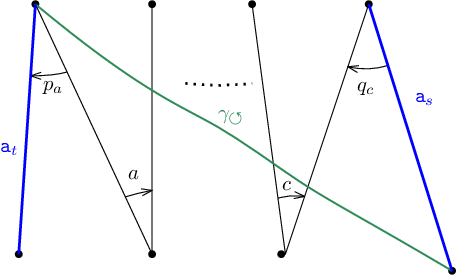}
    \hskip.5cm
\includegraphics[scale=0.6]{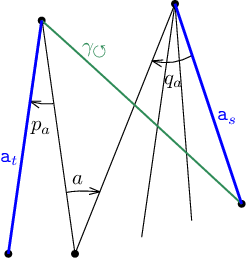}
\caption{Possible cases when $f_s(w_\gamma)$ is obtained by adding a hook on the right.}
    \label{fig:add-hook-on-right}
\end{figure}

{\bf Case b.} $w_\gamma$ is not direct and $f_s(\gamma)$ is obtained by removing a cohook $b q_b^{-1}$ on the right.

In this case, since $M$ is not injective, we must be able to add a hook $p_a^{-1} a$ on the left of $w_\gamma$. Figure~\ref{fig:hook-left-cohook-right} illustrates this case. It is clear that $\tau^{-1} (\gamma) = f_t(f_s(\gamma))$. 

\begin{figure}[ht!]
    \centering
\includegraphics[scale=0.7]{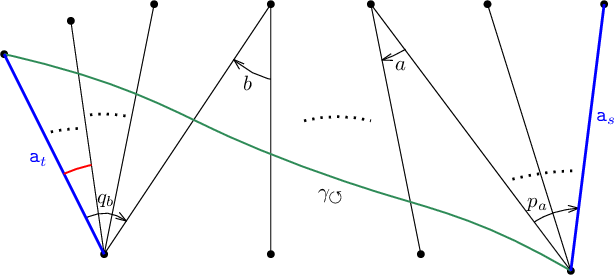}
    \caption{Possible case when $f_s(w_\gamma)$ is obtained by removing a cohook on the right and $M$ is not injective.}
    \label{fig:hook-left-cohook-right}
\end{figure}

This finishes the proof of $(2)$ and $(3)$. Condition $(4)$ is clear since $\tau^{-1}$ acts as the identity on band modules.
\end{proof}

\begin{figure}[ht!]
    \centering
    \includegraphics[scale=0.6]{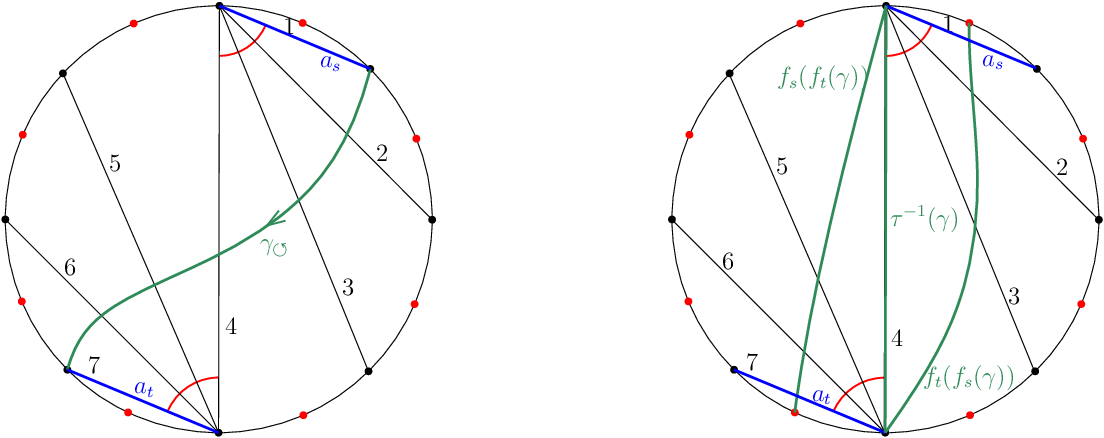}
    \caption{Example where $\tau^{-1} (\gamma), f_s(f_t(\gamma))$ and $f_t(f_s(\gamma))$ are pairwise non-equivalent.}
    \label{fig:not-commutative}
\end{figure}

\begin{remark}
    Note that $\tau^{-1} (\gamma)$ is not necessarily equivalent to $f_s (f_t (\gamma))$ or to $f_t(f_s (\gamma))$, as pointed out in Remark~\ref{rmk:not-commutative}. Figure~\ref{fig:not-commutative} shows permissible arcs corresponding to the strings $w, f_s(f_t(w))$ and $f_t(f_s(w))$, none of which are zero arcs. We have $f_s(f_t(w)) \neq f_t(f_s(w))$. Moreover, $\tau^{-1}(w)$ corresponds to the injective module at vertex $4$ and $\tau^{-1} (\gamma)$ is the arc in $\sP$ corresponding to this vertex. 
\end{remark}

\section{Classification of support tau-tilting modules }\label{sec:support-tau-tilting}

In this section we will use the geometric model of string algebras to give a combinatorial classification of support $\tau$-tilting modules over this class of algebras. We begin by giving the necessary background on $\tau$-tilting theory, and we refer the reader to~\cite{AIR} for more details. 

\begin{definition}\label{def:crossing-M}
    Let $A$ be a finite-dimensional algebra and $M \in \mod{A}$. Denote by $|M|$ the number of non-isomorphic indecomposable direct summands of $M$. 
    \begin{enumerate}
    \item $M$ is {\it $\tau$-rigid} if $\Hom_A (M, \tau M) = 0$.
    \item $M$ is {\it $\tau$-tilting} if $M$ is $\tau$-rigid and $|M|=|A|$.
    \item $M$ is {\it support $\tau$-tilting} if there is an idempotent $e \in A$ such that $M$ is a $\tau$-tilting $(A/ \left< e \right>)$-module.
    \item A {\it support $\tau$-rigid pair} is a pair $(P,M)$ consisting of a projective $A$-module $P$ and a $\tau$-rigid $A$-module $M$ such that $\Hom_A(P,M)=0$.
    \item A {\it support $\tau$-tilting pair} is a $\tau$-rigid pair $(P, M)$ for which $|M|+|P|=|A|$.
    \end{enumerate}
\end{definition}

\begin{proposition}\cite[Proposition 2.3]{AIR}
    Let $A$ be a finite-dimensional algebra and $M \in \mod{A}$. The module $M$ is support $\tau$-tilting if and only if there is a projective $A$-module $P$ such that $(P,M)$ is a support $\tau$-tilting pair. 
\end{proposition}

\begin{proposition}\cite[Corollary 2.13]{AIR}\label{prop:AIR-maximal}
    A pair $(P,M)$ is support $\tau$-tilting if and only if it is maximal with respect to the $\tau$-rigid pair property.
\end{proposition}

Note that if the algebra $A$ is a string algebra, then no band module is $\tau$-rigid since the Auslander-Reiten translate acts as the identity. Therefore, we will only need to consider string modules. 

\subsection{Crossings of $\sM$-arcs.} As we have seen in Section~\ref{sec:morphisms}, we can describe morphisms between two string modules $M_\beta$ and $M_\gamma$ in terms of the crossings and intersections between the corresponding $\sR$-arcs $\beta$ and $\gamma$. In this section we will see that certain representatives of $\sM$-arcs, namely the clockwise-most ones, are more useful to describe $\tau$-rigidity, in the sense that most crossings between these arcs detect morphisms in the Hom-space $\Hom(M,\tau N)$, where $M, N$ are string modules. 

In order to see this, we first need to define an orientation for each crossing. 

\begin{figure}
\centering
\includegraphics[scale=.6]{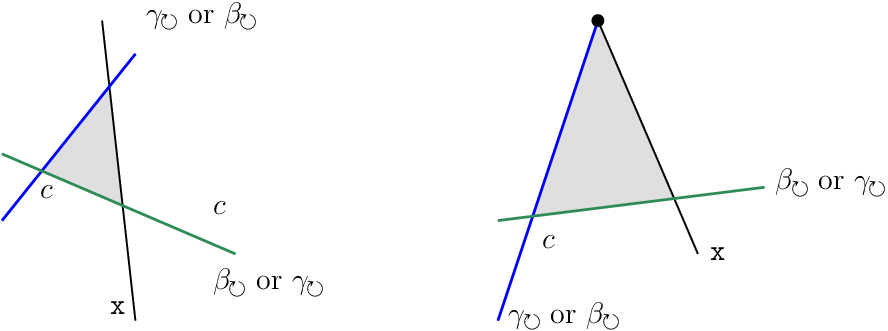}
\caption{Each crossing of permissible clockwise-most $\sM$-arcs defines a triangle.}
\label{fig:crossing-triangle}
\end{figure}

\begin{lemma}\label{lem:triangle}
    Let $\gamma$ and $\beta$ be non-trivial permissible clockwise-most $\sM$-arcs, and suppose they cross at a point $c$ in their interior. Then there is one arc $\arc{x}$ in $\sP$ and a triangle of one of the forms given in Figure~\ref{fig:crossing-triangle} with no puncture in $\sR$ in its interior. 
\end{lemma}
\begin{proof}
Let $\Delta$ be the tile the crossing $c$ lies in. If there is $\arc{x} \in \sP$ which bounds $\Delta$ which is intersected by $\gamma$ and $\beta$, then we clearly have a triangle as the one on the left of Figure~\ref{fig:crossing-triangle}. 
    
Suppose the segments of $\gamma$ and $\beta$ in $\Delta$ where the crossing $c$ occurs don't have a common intersection with $\sP$. Since they are permissible arcs which cross each other, at least one of the segments must be an end-segment. The fact that there is a local triangle $\delta$ follows from the fact that both arcs are permissible and clockwise-most $\sM$-arcs. Moreover, one of the vertices $p$ of $\delta$ must be in $\sM$. Let $\arc{x}$ be the arc in $\sP$ bounding $\delta$. If the end-segment of $\gamma$ or $\beta$ follows $\arc{x}$ around $p$ in the anticlockwise order instead, then we must be in one of the situations presented in Figure~\ref{fig:triangle-delta}. Therefore, there is another triangle $\delta'$ of the form given on the right of Figure~\ref{fig:crossing-triangle}.  

\begin{figure}[ht!]
\centering
\includegraphics[scale=0.6]{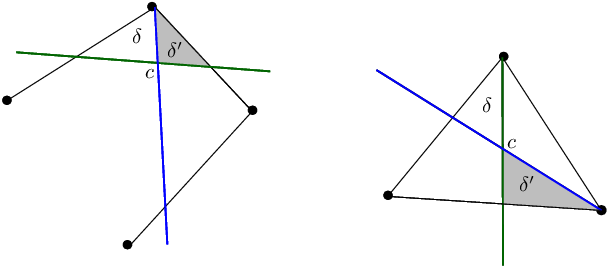}
\caption{The case when $\gamma$ and $\beta$ don't have a common intersection with $\sP$.}
\label{fig:triangle-delta}
\end{figure}

The fact that the red point of $\Delta$ does not lie inside the triangle follows from the definitions of permissible arc and clockwise-most $\sM$-arc. 
\end{proof}

The existence of this triangle permits us to define an orientation for each crossing.

\begin{definition}
    Let $\gamma$ and $\beta$ be non-trivial permissible clockwise-most $\sM$-arcs, and suppose they cross at a point $c$ in their interior. Let $\Delta$ be a triangle as in Figure~\ref{fig:crossing-triangle} associated to $c$. We say this crossing is from $\beta$ to $\gamma$ if $\gamma$ follows $\beta$ immediately in the clockwise order around the triangle. 
\end{definition}

\begin{remark}
    The orientation of the crossing is well defined. Indeed, if a crossing is incident with two triangles as in Figure~\ref{fig:crossing-triangle}, then they are either opposite to each other, in which case the orientation is the same for both triangles, or they are adjacent. Since $\sP$ is a dissection, it follows that both triangles share a marked point in $\sM$. But then, at least one of the triangles wouldn't be of the form given by the right picture in Figure~\ref{fig:crossing-triangle}, a contradiction.
\end{remark}

\subsection{Support $\tau$-tilting via crossings}
Let $M, N$ be string modules and $\gamma_M, \gamma_N$ be the corresponding permissible clockwise-most $\sM$-arcs. We will now give a formula for the dimension of the Hom-space $\Hom (M, \tau N)$ in terms of the crossings between $\gamma_M$ and $\gamma_N$. This has been done in the context of skewed-gentle algebras in \cite{HZZ}. We note that in \cite{HZZ}, the authors have the opposite orientation when defining the algebra from the surface. This implies that they use the anticlockwise-most arcs instead. 

\begin{lemma}\label{lem:R-arc-crossing-implies-clockwise-crossing}
Let $\beta$ and $\gamma$ be two permissible $\sR$-arcs and $\beta_\circlearrowright$ and $\gamma_\circlearrowright$ be the corresponding permissible clockwise-most $\sM$-arcs. If there is a crossing from $\beta$ to $\gamma$, then there is a crossing from $\beta_\circlearrowright$ to $\gamma_\circlearrowright$.
\end{lemma}
\begin{proof}
This statement follows immediately from Definition~\ref{def:crossing-R} and Definition~\ref{def:crossing-M}, since at least one of the segments where the crossing between $\beta$ and $\gamma$ occurs is an interior segment. Note that the triangle associated to the crossing from $\beta_\circlearrowright$ to $\gamma_\circlearrowright$ is as the one on the left in Figure~\ref{fig:crossing-triangle}.
\end{proof}

\begin{lemma}\label{lem:morphism-implies-crossing}
    Let $M$ and $N$ be string $A$-modules, and denote by $\gamma_M$ and $\gamma_N$ the corresponding permissible clockwise-most $\sM$-arcs. If $\Hom_A (M, \tau N) \neq 0$ then there is a crossing from $\gamma_M$ to $\gamma_N$. 
\end{lemma}
\begin{proof}
    Since $\Hom (M, \tau N) \neq 0$, we have in particular that $N= \tau^{-1} (\tau N) \neq 0$ and $\tau N$ is not injective. Moreover, $w_M$ and $w_{\tau N}$ are of the form: 
    \[
    w_M = (w_M)_\ell a^{-1} e b (w_M)_r, \text{ and } w_{\tau N} = (w_{\tau N})_\ell c e d^{-1} (w_{\tau N})_r,
    \]
    with $a, b, c, d \in Q_1$ if $(w_M)_\ell a^{-1}, b (w_M)_r, (w_{\tau N})_\ell c, d^{-1}(w_{\tau N})_r \neq 0$, respectively. 

    In order to check whether $\gamma_M$ crosses $\gamma_N$, we need to compute $\tau^{-1} (\tau N)$. 

    {\bf Case 1:} $f_t(w_{\tau N}) = 0$. 

Then $w_{\tau N}$ must be a left maximal direct path. Since $\tau N$ is not injective, there must be an arrow $d'$ such that $w_{\tau N} d'^{-1}$ is a string and so $f_s (w_{\tau N}) = w_{\tau N} d'^{-1}p_{d'}$, where $d'^{-1}p_{d'}$ is the hook of $d'$. It then follows that $w_N = f_t f_s (w_{\tau N}) = p_{d'}$.

Since $w_{\tau N}$ is a left maximal direct path, we must have $d^{-1} (w_{\tau N})_r = 0$ and $(w_{\tau N})_\ell c e = q_{d'}$, i.e. $(w_{\tau N})_\ell c e d'^{-1}$ is the cohook of $d'$. Figure~\ref{fig:ft-is-0} illustrates this situation. In particular, we have indeed a crossing from $\gamma_M$ to $\gamma_N$. 
    
\begin{figure}[ht!]
\centering
\includegraphics[scale=.5]{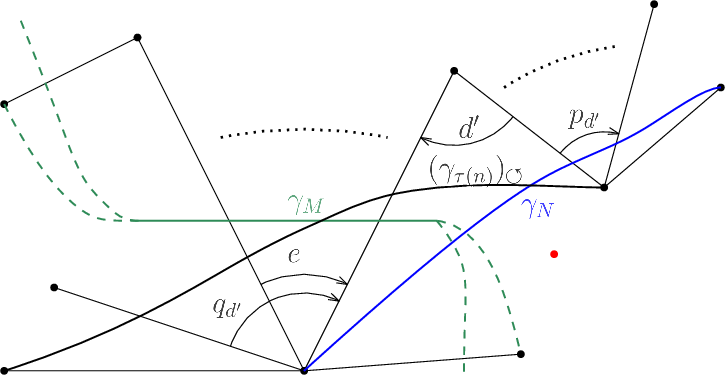}
\caption{The case where $f_t(w_{\tau N}) = 0$.}
\label{fig:ft-is-0}
\end{figure}

{\bf Case 2:} $f_s(w_{\tau N}) = 0$. This case is similar to the previous case.

{\bf Case 3:} $f_t(w_{\tau N}) \neq 0$ and $f_s(w_{\tau N}) \neq 0$. 

In this case we have $w_N = f_s(f_t (w_{\tau_N})) = f_t(f_s (w_{\tau_N}))$. We will calculate the former. 

{\it Case 3.1:} $w_{\tau N}$ is not left maximal or $(w_{\tau N})_\ell$ has a letter in $Q_1^{-1}$. 

Then $f_t(w_{\tau N}) = w'c e d^{-1} (w_{\tau N})_r$, for some non-zero string $w'$.

{\it Case 3.1.1:} $w_{\tau N}$ is not right maximal or $(w_{\tau N})_r$ has a letter in $Q_1$. 

Then $w_N = f_s(f_t(w_{\tau N})) = w' c e d^{-1} w''$, for some non-zero string $w''$. In particular, there is a non-zero morphism from $M$ to $N$ associated to the overlap $e$, and $e$ is neither a prefix nor a suffix of $w_N$. Hence, there is a crossing from the $\sR$-arc of $M$ to the $\sR$-arc of $N$. The result then follows from Lemma~\ref{lem:R-arc-crossing-implies-clockwise-crossing}. 

{\it Case 3.1.2:} $(w_{\tau N})_r$ is a right maximal inverse string.

Suppose $e$ has a letter in $Q_1$ and let $x$ denote the last such letter, i.e.~$ed^{-1}(w_{\tau N})_r = e'x q_x^{-1}$, where $x q_x^{-1}$ is the cohook of $x$. Then $w_N = f_s (f_t (w_{\tau N})) = w'c e'$. 

Therefore, we have a non-zero morphism from $M$ to $N$ associated to the overlap $e'$. Moreover, $e'$ is not a suffix of $w_M$ nor a prefix of $w_N$. Hence, there must be a crossing from the $\sR$-arc of $M$ to the $\sR$-arc of $N$, and the result follows again from Lemma~\ref{lem:R-arc-crossing-implies-clockwise-crossing}. 

So suppose $e$ is an inverse string. Then $ced^{-1}(w_{\tau N})_r$ is the cohook of $c$ and $w_N = f_s(f_t(w_{\tau N})) = w'$. We have no morphism from $M$ to $N$ associated to $e$ or a subpath of $e$, but there is a crossing from $\gamma_M$ to $\gamma_N$, as illustrated in Figure~\ref{fig:e-inverse}. 

\begin{figure}[ht!]
\centering
\includegraphics[scale=.6]{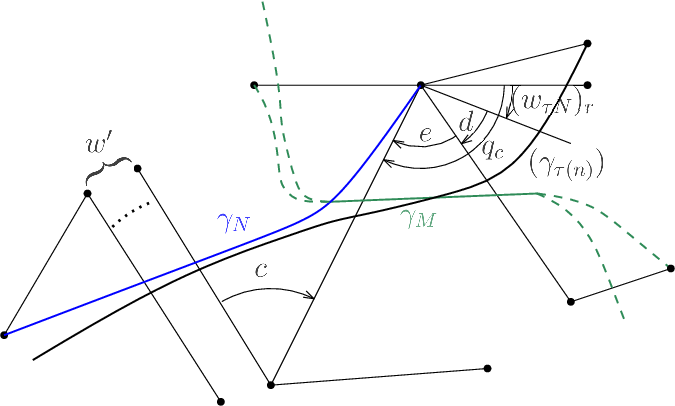}
\caption{The case where the cohook of $c$ is removed on the right.}
\label{fig:e-inverse}
\end{figure}

{\it Case 3.2:} $(w_{\tau N})_\ell$ is a left maximal direct string. 

{\it Case 3.2.1:} $e$ has a letter in $Q_1^{-1}$. Denote by $x^{-1}$ the first such letter. Then $(w_{\tau N})_\ell e=q_x x^{-1} e'$, where $x q_x^{-1}$ is the cohook of $x$. Thus $f_t(w_{\tau N}) = e'd^{-1} (w_{\tau N})_r$. 

{\it Case 3.2.1.1:} $(w_{\tau N})_r$ is not right maximal or if it has a letter in $Q_1$.

Then $w_N = f_s(f_t(w_{\tau N})) = e'd^{-1}w''$, for some non-zero string $w''$. Hence, there is a non-zero morphism from $M$ to $N$ associated to the overlap $e'$. Moreover, $e'$ is not a suffix of $w_N$ nor a prefix of $w_M$, meaning that there is a crossing from the $\sR$-arc of $M$ to the $\sR$-arc of $N$ and the result follows from Lemma~\ref{lem:R-arc-crossing-implies-clockwise-crossing}. 

{\it Case 3.2.1.2:} $(w_{\tau N})_r$ is a right maximal inverse string. 

Since $f_s(f_t(w_{\tau N})) = w_N \neq 0$, $e'$ must have a letter in $Q_1$. Let $y$ be the last such letter, and write $e' = e'' y q^{-1}$, where $q$ is a direct path and $e''$ is a non-zero string. Then $y q^{-1} d^{-1} (w_{\tau N})_r$ is the cohook of $y$ and $w_N = e''$. Again, there is a non-zero morphism from $M$ to $N$ associated to the overlap $e''$, which is neither a prefix  nor a suffix of $w_M$, and so the result follows from Lemma~\ref{lem:R-arc-crossing-implies-clockwise-crossing}. 

{\it Case 3.2.2:} $e$ is a direct path. 

Since $f_t(w_{\tau N}) \neq 0$, we must have that $d^{-1} (w_{\tau N})_r$ is non-zero and so $d ((w_{\tau N})_\ell c e)^{-1}$ is the cohook of $d$ and $f_t(w_{\tau N}) = (w_{\tau N})_r$. Given that $f_s(f_t(w_{\tau N})) \neq 0$, we can either add a hook on the right of $(w_{\tau N})_r$ or $(w_{\tau N})_r$ has a letter in $Q_1$. This case is similar to the case illustrated in Figure~\ref{fig:e-inverse}, and the crossing from $\gamma_M$ to $\gamma_N$ takes place in the tile associated to the arrow $d$.
\end{proof}

\begin{figure}
\centering
\includegraphics[scale=.6]{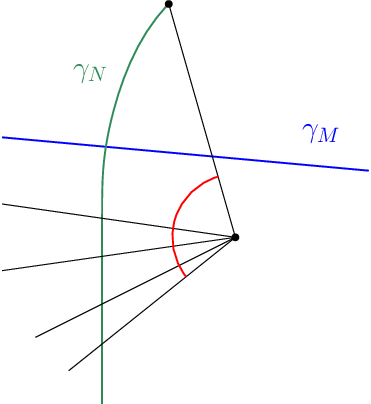}
\caption{Crossings from $\gamma_M$ to $\gamma_N$ not associated with morphisms from $M$ to $\tau N$.}
\label{fig:crossing-X}
\end{figure}

\begin{definition}\label{good-crossing}
Let $M, N$ be two string modules and $\gamma_M, \gamma_N$ be the corresponding clockwise-most $\sM$-arcs. A crossing from $\gamma_M$ to $\gamma_N$ which is, up to isotopy, of the form shown in Figure~\ref{fig:crossing-X} is called a {\it good crossing}. Given such a crossing, we have the following conditions:
\begin{itemize}
    \item $w_N = w_\ell a_1 \ldots a_{r-1}$, where $w_\ell$ is a non-zero string,  and $a_1 \ldots a_{r-1}$ is a direct path of length $\geq 1$;
    \item there exists $a_r \in Q_1$ such that $a_1 \ldots a_{r-1} a_r$ corresponds to a label;
    \item $\gamma_M$ intersects the arc in $\sP$ corresponding to $t(a_r)$.
\end{itemize}
We denote by $\sX_{M,N}$ the number of good crossings from $\gamma_M$ to $\gamma_N$. 
\end{definition}

The following proposition gives the dimension of the Hom-space $\Hom (M, \tau N)$, where $M$ and $N$ are string modules, in terms of crossings of clockwise-most $\sM$-arcs. The case in common with~\cite{HZZ} is when the algebra is gentle. In this case, there are no good crossings and the dimension of $\Hom (M, \tau N)$ is precisely the number of crossings from $\gamma_M$ to $\gamma_N$. 

\begin{proposition}\label{prop:crossings-tau-extensions}
    Let $M, N$ be string $A$-modules, and denote by $\gamma_M$ and $\gamma_N$ the corresponding permissible clockwise-most $\sM$-arcs. We have
    \[
    \dim \Hom_A (M,\tau N) = \sC (\gamma_M, \gamma_N) -\sX_{M,N}.
    \]
\end{proposition}
\begin{proof}
    Take a crossing $c$ from $\gamma_M$ to $\gamma_N$ which is not a good crossing. We will first associate a non-zero morphism from $M$ to $\tau N$ to the crossing $c$.
    
    If $c$ defines a triangle as the one on the left in Figure~\ref{fig:crossing-triangle}, then $\gamma_M$ and $\gamma_N$ have a common intersection with $\sP$ and  we have $w_M = (w_M)_\ell a^{-1} e b (w_M)_r$ and $w_N =  (w_N)_\ell c e d^{-1} (w_N)_r$, where the overlap $e$ is defined by the crossing $c$ as illustrated in Figure~\ref{fig:e-corresponding-to-crossing} and where $a, b, c, d \in Q_1$ are such that at least one of $(w_M)_\ell a^{-1}$ and $(w_N)_\ell c$ (resp. one of $b (w_M)_r$ and $d^{-1} (w_N)_r$) is non-zero. 

    In order to calculate $w_{\tau N}$, we  add a cohook or otherwise remove a hook, both at the start and end of $w_N$. 

    {\bf Case 1:} $(w_N)_\ell c \neq 0$ and $d^{-1} (w_N)_r \neq 0$.

    In this case, $e$ is still a substring of $w_{\tau N}$. Therefore there is a non-zero morphism from $M$ to $\tau N$ associated to the overlap $e$. 

    {\bf Case 2:} $(w_N)_\ell c \neq 0$ but $d^{-1} (w_N)_r = 0$.

    Then we have $b (w_M)_r \neq 0$. Moreover, since $c$ is not a good crossing, $f_s(w_N)$ is obtained by adding the cohook of $b$ at the end of $w_N$, i.e.~$f_s(w_N) = w_N b q_b^{-1}$. Since $c \neq 0$, it follows that there is a morphism from $M$ to $\tau N$ corresponding to the overlap $ebp^{-1}$, where $p^{-1}$ is a prefix of $q_b^{-1}$. 

    {\bf Case 3:} $(w_N)_\ell c = 0$ and $d^{-1} (w_N)_r \neq 0$. 

    This case is similar to Case 2. 

    {\bf Case 4:} $(w_N)_\ell c = 0$ and $d^{-1} (w_N)_r = 0$. 

    Then we have $a, b \neq 0$ and $a^{-1} e b$ is a string. Therefore $w_{\tau N} = q_a a^{-1} e b q_b^{-1}$ and we have a non-zero morphism from $M$ to $\tau N$ with associated overlap $p_1a^{-1} e b p_2^{-1}$, where $p_1$ is a suffix of $q_a$ and $p_2^{-1}$ is a prefix of $q_b^{-1}$. 

  Now suppose $c$ defines a triangle as the one on the right in Figure~\ref{fig:crossing-triangle}. The only case which is not isotopic to the crossings analysed above is presented in Figure~\ref{fig:crossing-not-common-P}.
  
  \begin{figure}[ht!]
      \centering
      \includegraphics[scale=0.8]{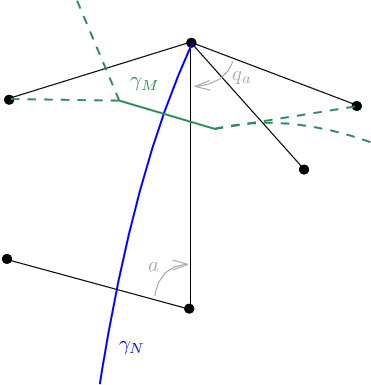}
      \caption{A crossing from $\gamma_M$ to $\gamma_N$ not defining a triangle as the one on the left in Figure~\ref{fig:crossing-triangle}.}
      \label{fig:crossing-not-common-P}
  \end{figure}
  
  In this case, we have an arrow $a \in Q_1$ such that $s(a) = t(w_N)$ and $w_N a$ is a string, since the $c$ is not a good crossing. Therefore $f_s(w_N)$ is obtained by adding the cohook $aq_a^{-1}$ on the right of $w_N$, and $q_a^{-1}$ a substring of $w_{\tau N}$. Hence, there is a non-zero morphism from $M$ to $\tau N$ associated to an overlap which is a prefix of $q_a^{-1}$.
  
  Now let $e$ be an overlap giving rise to a non-zero morphism from $M$ to $\tau N$. By Lemma~\ref{lem:morphism-implies-crossing}, there is a crossing from $M$ to $N$. If it is a good crossing then one of the segments must be an end-segment. The only cases in the proof of Lemma~\ref{lem:morphism-implies-crossing} where this occurs are the cases 1, 2, 3.1.2 and 3.2.2. But as we can see from Figures~\ref{fig:ft-is-0} and~\ref{fig:e-inverse}, which illustrate these cases, the corresponding crossing is not of form given in Figure~\ref{fig:crossing-X}. 
  
  It is clear that the mappings defined above between basis elements of $\Hom (M, \tau N)$ and crossings from $\gamma_M$ to $\gamma_N$ which are not good are inverse of each other, and so the result follows.
\end{proof}

The following corollary follows immediately from Proposition~\ref{prop:crossings-tau-extensions} since each crossing has a unique orientation.

\begin{corollary}\label{cor:crossings-tau-extensions}
    Let $M, N$ be string $A$-modules, and denote by $\gamma_M$ and $\gamma_N$ the corresponding permissible clockwise-most $\sM$-arcs. We have
    \[
    \dim \Hom_A (M,\tau N) +\dim \Hom_A (N, \tau M) = \sC (\gamma_M, \gamma_N) +\sC (\gamma_N, \gamma_M) -\sX_{M,N} - \sX_{N,M}.
    \]
\end{corollary}

\begin{definition}
   A {\it generalised permissible arc} is either an arc in $\sP$ or a permissible clockwise-most $\sM$-arc.  
\end{definition}

We are now ready to give a combinatorial classification of support $\tau$-tilting modules over a string algebra, which recovers the case where the algebra is gentle (cf.~\cite{HZZ}). We note that there is a different classification of support $\tau$-tilting modules over string algebras in \cite{EJR}.

\begin{theorem}\label{thm:tau-tilt}
Let $A$ be a string algebra. There is a one-to-one correspondence between the set of support $\tau$-tilting pairs and the set of maximal collections of generalised permissible arcs for which each crossing is good.
\end{theorem}
\begin{proof}
This statement follows immediately from Proposition~\ref{prop:crossings-tau-extensions} and Proposition~\ref{prop:AIR-maximal}. The arcs in $\sP$ correspond to the summands of the projective module in the pair and the remaining arcs correspond to the summands of the $\tau$-rigid module in the pair. 
\end{proof}

\begin{figure}[ht!]
    \centering
    \includegraphics[scale=0.7]{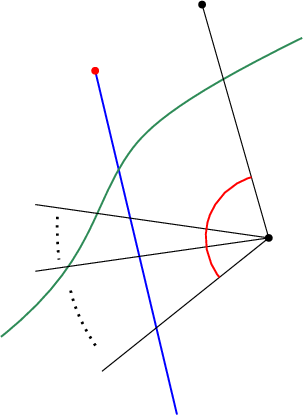}
    \caption{Crossings between $\sR$-arcs allowed in support $\tau$-tilting modules.}
    \label{fig:crossing-X-red}
\end{figure}

\begin{figure}[ht!]
\centering
\includegraphics[scale=.7]{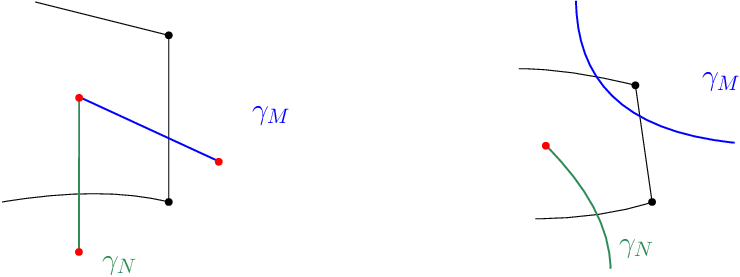}
\caption{Two local non-crossing configurations which are not allowed unless we are in the situation of Figure~\ref{fig:non-crossing-configuration-label}.}
    \label{fig:non-crossing-configuration}
\end{figure}

\begin{figure}[ht!]
\centering
\includegraphics[scale=.7]{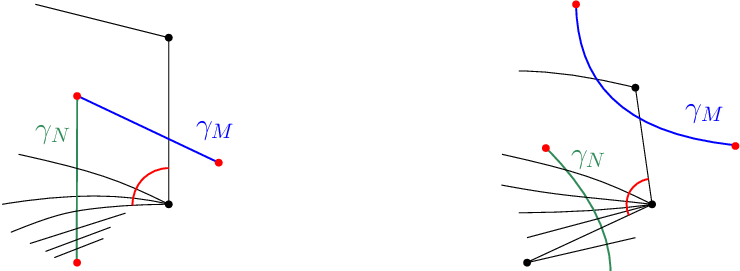}
\caption{A local configuration which is allowed.}
    \label{fig:non-crossing-configuration-label}
\end{figure}

\begin{remark}
    We can reformulate this classification in terms of $\sR$-arcs as follows. A collection of arcs in $\sP$ together with permissible $\sR$-arcs corresponds to a support $\tau$-tilting pair if and only if it is maximal with respect to the following properties:
    \begin{enumerate}
        \item all crossings are of the form given in Figure~\ref{fig:crossing-X-red}, and
        \item if two $\sR$-arcs satisfy the condition shown in Figure~\ref{fig:non-crossing-configuration}, then there is a label as in Figure~\ref{fig:non-crossing-configuration-label}. 
    \end{enumerate}
\end{remark}

\begin{remark}
Consider the algebra $G_1$ given by the quiver $\xymatrix{1 \ar[r]^a & 2 \ar[r]^b & 3}$, bound by the relation $ab$. Then $G_1$ is gentle and it is the tiling algebra of the surface in Figure~\ref{fig:tilted-A}. 

\begin{figure}
\centering
\includegraphics[scale=0.6]{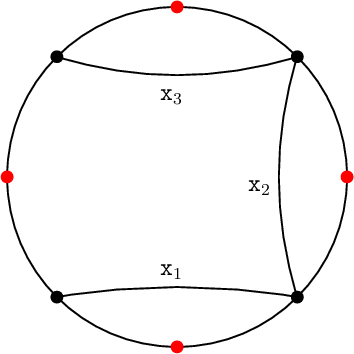}
\caption{Surface with tiling algebra $G_1$.}
\label{fig:tilted-A}
\end{figure}

Consider the projective module $P_2$ and the simple module $S_1$. We have that $P_2 \oplus S_1$ is not $\tau$-rigid, since $\Hom_{G_1} (P_2, \tau S_1) \neq 0$. In fact, the dimension of $\Hom_{G_1} (P_2, \tau S_1)$ is one. 

On the one hand, we have no crossings or intersections between the permissible $\sR$-arcs associated to $P_2$ and $S_1$, and on the other hand, as expected from Proposition~\ref{prop:crossings-tau-extensions}, there is one crossing between the corresponding clockwise-most permissible $\sM$-arcs; see Figure~\ref{fig:tilted-R-M}. 

\begin{figure}[ht!]
\centering
\includegraphics[scale=0.6]{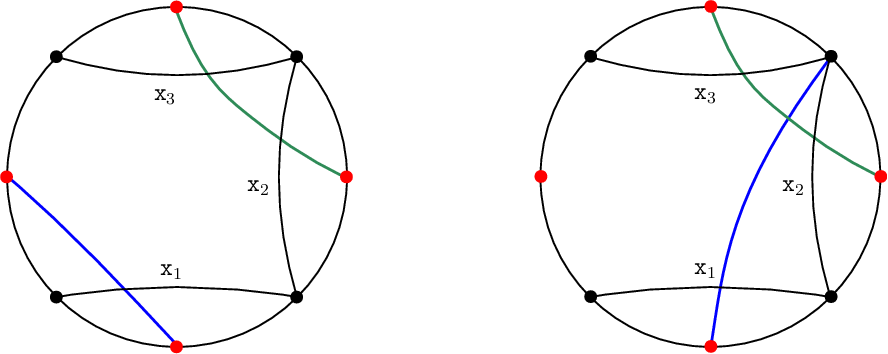}
\caption{The $\sR$-arcs and clockwise-most $\sM$-arcs associated to $S_1$ and $P_2$.}
\label{fig:tilted-R-M}
\end{figure}

Now consider the algebra $G_2$ given by the quiver $\xymatrix{1 \ar[r]^a & 2 \ar[r]^b & 3 \ar[r]^c & 4}$ bound by the relation $abc$. Figure~\ref{fig:G2-example} shows the $\sR$-arcs associated to the projective modules $P_1$ and $P_2$. We clearly have that $P_1 \oplus P_2$ is a $\tau$-rigid $G_2$-module.

\begin{figure}[ht!]
    \centering
    \includegraphics[scale=0.6]{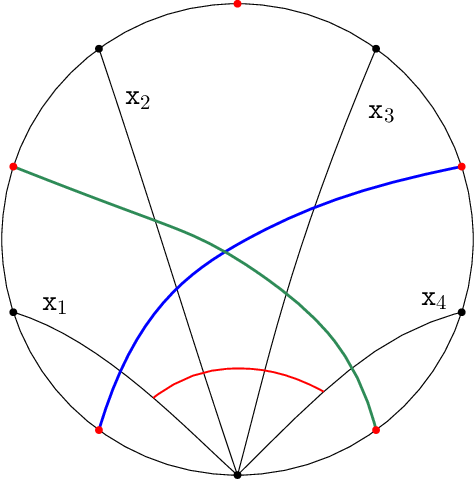}
    \caption{The surface associated to $G_2$ and the $\sR$-arcs corresponding to $ab$ and $bc$.}
    \label{fig:G2-example}
\end{figure}

These examples illustrate the reason why we consider $\sM$-arcs to describe the Auslander-Reiten translate and classify support $\tau$-tilting modules. On the one hand, $\sM$-arcs allow us to avoid considering local configurations which are not given by crossings or intersections. On the other hand, using $\sR$-arcs in the description of $\tau$-rigidity does not exclude certain crossings.
\end{remark}

\begin{example}
   Consider the algebra given in Example~\ref{ex:labeled}. In Figure~\ref{fig:ex-sec-5}, we have the collection of clockwise-most $\sM$-arcs together with the highlighted arc in $\sP$ associated to the support $\tau$-tilting pair $(P_3, S_1 \oplus M(a_5) \oplus M(a_5^{-1}a_6) \oplus M(a_1a_2))$. Note that there is a self-crossing of the arc associated to $M(a_1a_2)$. However, this crossing is of the form shown in Figure~\ref{fig:crossing-X}, and indeed $M(a_1a_2) = P_1$ is a $\tau$-rigid module.

\begin{figure}[ht!]
\centering 
\includegraphics[scale=0.8]{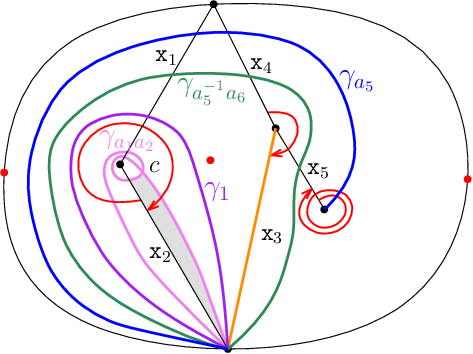}
\caption{An example of a support $\tau$-tilting pair.}
\label{fig:ex-sec-5}
\end{figure}
\end{example}


\bibliographystyle{plain}
\bibliography{biblio}

\end{document}